\documentclass{amsart}
\usepackage{amssymb}
\usepackage[arrow,matrix]{xy}

\newcommand\TR{\operatorname{T}}

\newcommand\End{\operatorname{End}}
\newcommand\END{\operatorname{END}}
\newcommand\ENDG{\END_{\GR}(E)}
\newcommand\Hom{\operatorname{Hom}}
\newcommand\CC{\mathbb C}

\newcommand\RR{\mathbb R}
\newcommand\ZZ{\mathbb Z}
\newcommand\pa{\partial}
\newcommand\supp{\operatorname{supp}}

\newcommand\CI{{\mathcal C}^{\infty}}
\newcommand\CIc{{\mathcal C}^{\infty}_{\text{c}}}
\newcommand\ie{i.e.\/ }
\newcommand\hden{{\Omega^{\lambda}_d}}
\newcommand\VD{{\mathcal D}}

\newcommand{\Gr}[1]{{\mathcal G}^{(#1)}}
\newcommand{\GR}{\mathcal G}
\newcommand{\LGR}{\mathcal L}

\newcommand{\PS}[1]{\Psi^{#1}(\GR;E)}
\newcommand{\tPS}[1]{\Psi^{#1}(\GR)}
\newcommand{\ttPS}[1]{\Psi_{\loc}^{#1}(\GR;E)}
\newcommand{\tttPS}[1]{\Psi_{\loc}^{#1}(\GR)}

\newcommand{\FAM}{P=(P_x,x \in \Gr0)}

\newcommand{\loc}{\operatorname{loc}}
\newcommand{\cl}{\operatorname{cl}}
\newcommand{\A}{s}
\newcommand{\prop}{\operatorname{prop}}
\newcommand{\comp}{\operatorname{comp}}
\newcommand{\adb}{\operatorname{adb}}
\newcommand{\dist}{\operatorname{dist}}

\newcommand{\alp }{r }
\newcommand{\bet }{d }
\newcommand{\gm }{\Gamma }
\newcommand{\lon }{\longrightarrow }
\newcommand{\be }{\begin{eqnarray*}}
\newcommand{\ee }{\end{eqnarray*}}
\newcommand{\GGR}{{\GR}}

\def\ie{i.e.\ }

\newcommand{\frakg}{{\mathfrak g}}

\newtheorem{theorem}{Theorem}
\newtheorem{proposition}{Proposition}
\newtheorem{corollary}{Corollary}
\newtheorem{lemma}{Lemma}
\newtheorem{definition}{Definition}

\theoremstyle{remark}

\newtheorem{example}{Example}
\begin{document}

\title[Pseudodifferential operators on groupoids]
{Pseudodifferential operators on differential groupoids}

\author[V. Nistor]{Victor Nistor}
\address{Department of Mathematics, Pennsylvania State University}
\email{nistor@math.psu.edu}
\thanks{
Nistor is partially supported by NSF Young Investigator Award
DMS-9457859 and a Sloan 
research fellowship, preprints available from
{\bf http:{\scriptsize//}www.math.psu.edu{\scriptsize/}nistor{\scriptsize/}}
}

\author[A. Weinstein]{Alan Weinstein}
\address{Department of Mathematics, University of California, Berkeley}
\email{alanw@math.berkeley.edu}
\thanks{Weinstein is partially supported by NSF Grants DMS-9309653 and DMS-9625122}

\author[Ping Xu]{Ping Xu}
\address{Department of Mathematics, Pennsylvania State University}
\email{ping@math.psu.edu}
\thanks{Xu is partially supported by NSF Grant DMS 95-04913}

\date{\today} 

\begin{abstract} 
We construct an algebra of pseudodifferential operators on each groupoid 
in a class that generalizes differentiable groupoids to allow manifolds 
with corners. We show that this 
construction encompasses many examples. The subalgebra of 
regularizing operators is identified with the smooth algebra of the 
groupoid, in the sense of non-commutative geometry. Symbol calculus
for our algebra lies in the Poisson algebra of functions on the dual
of the Lie algebroid of the groupoid. As applications, we give a new
proof of the Poincar\'e-Birkhoff-Witt theorem for Lie algebroids and
a concrete quantization of the Lie-Poisson
structure on the dual $A^*$ of a Lie algebroid.
\end{abstract}

\maketitle
\tableofcontents

\section*{Introduction}

Certain important 
applications of pseudodifferential operators
require variants of the original definition.
Among the many examples one can find in the literature are
regular or adiabatic families of pseudodifferential operators
\cite{Atiyah-Singer4,Witten1} and
pseudodifferential operators along the leaves of foliations
\cite{Connes3,Connes1,Moore-Schochet1,Nistor4},
on coverings \cite{Connes-Moscovici2,Nistor2}
or on certain singular spaces 
\cite{Melrose42,Melrose46,Melrose-Nistor2,Melrose-Nistor1}.

Since these classes of operators share many common features, it is
natural to ask whether they can be treated in a unified way.  
In this paper we shall suggest  an answer to this question.
For any ``almost differential'' groupoid (a class  which
allows manifolds with corners), we construct
an algebra of pseudodifferential operators.
We then show that our construction recovers (almost) all the
classes described above (for operators on manifolds with
boundary our algebra is slightly smaller than the one defined in 
\cite{Melrose42}). We expect our results
to have applications to analysis on singular spaces, not only
manifolds with corners.  

Our construction and results owe a great deal 
to the previous work of several authors, especially
Connes \cite{Connes3} and Melrose \cite{Melrose42,Melrose45,Melroseb}.  
A hint of the direction we take was given at the end of 
\cite{Weinstein1}.  The basic idea of our construction is to
consider families of pseudodifferential operators along the fibers of
the domain (or source) map of the groupoid. 
More precisely, for any almost differentiable groupoid (see
Definition \ref{Almost.differentiable}) we consider the fibers $\GR_x=d^{-1}(x)$
of the domain map $d$ consisting of all arrows with domain $x$. It follows
from the definition that these fibers are 
smooth manifolds (without corners). The calculus of 
pseudodifferential operators on smooth manifolds is well understood and
by now a classical subject, see for example \cite{Hormander3}. 
We shall consider differentiable families of pseudodifferential operators 
$P_x$ on  the smooth manifolds $\GR_x$. Right translation by $g \in \GR$
defines an isomorphism $\GR_x \equiv \GR_y$ where $x$ is the domain
of $g$ and $y$ is the range of $g$. We say that the family $P_x$
is invariant if $P_x$ transforms to $P_y$ under the diffeomorphisms 
above (for all $g$). The algebra $\tPS{\infty}$
of {\em pseudodifferential operators on} $\GR$  that we shall consider
will consist of invariant differentiable families of operators $P_x$
as explained above (the actual definition also involves a technical
condition on the support of these operators). See Definition
\ref{Main.definition} for details. The relation with the work of 
Melrose relies on an alternative description of our algebra as an
algebra of distributions on $\GR$ with suitable properties
(compactly supported, and conormal with singular support contained in the 
set of units). This is contained in Theorem
\ref{Theorem.Full.Isom}. The difference between our theory and
Melrose's
lies in the fact that he
 considers a compactification of $\GR$ as a manifold with
corners, and his distributions are allowed to extend to the 
compactification, with precise behavior at the boundary. 
This is useful for the analysis of these operators. In contrast, our
work is purely algebraic (or geometric, depending on whether
one considers Lie algebroids as part of geometry or algebra). 

We now review the contents of the sections of this paper. In the first section
we recall the definitions of a groupoid, Lie algebroid and, the less
known definition of a {\em local} Lie groupoid. We extend the definition
of a Lie groupoid to include manifolds with corners. These groupoids 
are called {\em almost differentiable} groupoids. The second section
contains the definition of a pseudodifferential operator on a 
groupoid (really a family of pseudodifferential operators, as 
explained above) and the proof that they form an algebra, if a support
condition is included. We also extend this definition to include
local Lie groupoids. This is useful in the third section where we 
use this to give a new proof of the Poincar\'e-Birkhoff-Witt theorem for Lie
algebroids. In the process of proving this theorem we also exemplify
our definition of pseudodifferential operators on an almost differentiable
groupoid by describing the differential operators in this class.
As an application we give an explicit construction of a deformation
quantization of the Lie-Poisson structure on $A^*$, the dual of
Lie algebroid $A$.
The section entitled ``Examples'' contains just what the title 
suggests: for many particular examples of groupoids $\GR$ we explicitly
describe the algebra $\tPS{\infty}$ of pseudodifferential operators
on $\GR$. This recovers classes of operators that were previously 
defined using  {\em ad hoc} constructions. Our definition is  often
not only more general, but also simpler. This is the case for operators
along the leaves of foliations \cite{Connes1,Moore-Schochet1} or 
adiabatic families of operators. 

Since one of our main themes is that the Lie algebras of vector fields
which are central in \cite{Melrose-scatering} are in fact the spaces
of sections of Lie algebroids, 
we describe these Lie algebroids explicitly in each of our examples. 

In the sixth section of the paper, we describe the convolution kernels
(called {\em reduced} kernels) of operators in $\tPS{\infty}$. Then we
extend to our setting some fundamental results on principal symbols,
by reducing to the classical results. This makes our proofs short (and
easy). Finally, the last section treats the action of $\tPS{\infty}$
on functions on the units of $\GR$, and a few related topics.

The first author would like to thank Richard Melrose for
several useful conversations.

\section{Preliminaries}

In the following we allow  manifolds to have corners. 
Thus by ``manifold'' we shall mean a $C^\infty$ manifold, possibly 
with corners, and  by a ``smooth manifold'' we shall mean a manifold without
corners. By definition, if $M$ is a manifold with corners then
every point $p\in M$ has coordinate neighborhoods
diffeomorphic to $[0,\infty)^k \times \RR^{n-k}$. The transition functions
between such
coordinate neighborhoods must be smooth everywhere (including on the boundary). 
We shall use the following definition of submersions
between manifolds (with corners).

\begin{definition} A submersion between two manifolds with corners $M$ 
and $N$ is a differentiable map $f:M \to N$ such that 
$df_x:T_{x}M \to T_{f(x)}N$
is onto for any $x \in M$, and if $df_x(v)$ is an inward
pointing tangent vector to $N$, then $v$ is an inward pointing tangent 
vector to $M$. 
\end{definition}

The reason for introducing the definition above is that for any submersion  
$f: M \to N$, the set $M_y=f^{-1}(y)$, $y\in N$ is a {\em smooth} manifold,
just as for submersions of smooth manifolds. 

We shall study groupoids endowed with various
structures. (\cite{Renault1} is a general reference for some of what
follows.)  We recall first that a small
category is a category  whose class of morphisms is a set.
The class of objects of a small category is then a set as well.

\begin{definition} A groupoid  is a small category $\GR$ in which every 
morphism is invertible. 
\end{definition}

This is the shortest but least explicit definition.
We are going to make this definition more explicit in cases of interest.
The set of objects, or units, of $\GR$ will be denoted by 
$$
M=\GR^{(0)}=\operatorname{Ob}(\GR).
$$ 
The set of morphisms, or arrows, of $\GR$
will be denoted by 
$$
\GR^{(1)}=\operatorname{Mor}(\GR).
$$
We shall sometimes  
write $\GR$ instead of $\GR^{(1)}$ by abuse of notation.  For example,
when we consider a space of functions on $\GR$, we actually
mean a space of functions on $\GR^{(1)}$. 
We will denote by $d(g)$ [respectively
$r(g)$] the {\em domain} [respectively, the {\em range}] of the morphism 
$g:d(g) \to r(g).$
We thus obtain functions 
\begin{equation}
d,r:  \GR^{(1)} \longrightarrow  \GR^{(0)}
\end{equation}
that will play an important role bellow.
The  multiplication operator $\mu: (g, h) \mapsto\mu(g,h)=g h$
is defined on the set of composable pairs of
arrows $\GR^{(2)}$:
\begin{equation}\label{Mu.everywhere}
\mu:\GR^{(2)}=\GR^{(1)}\times_M \GR^{(1)}:=\{(g,h):d(g)=r(h)\} \longrightarrow
\Gr1.
\end{equation}
The inversion operation is a bijection 
$\iota:g\mapsto g^{-1}$ of $\Gr1$.
Denoting by $u(x)$ the identity morphism of the object  $x \in M= \Gr0$,
we obtain an inclusion of $\Gr0$ into $\Gr1$.  
We see that a  groupoid $\GR$ is completely determined by 
the spaces $\Gr0$ and $\Gr1$ and the structural morphisms $d,r,\mu,u,\iota$.
We sometimes write $\GR=(\Gr0,\Gr1,d,r,\mu,u,\iota)$. The structural
maps satisfy the following properties:
\\
(i)  $r(gh)=r(g)$, $d(gh)=d(h)$ for any pair $(g,h)\in \GR^{(2)}$,
and the partially defined multiplication $\mu$ is associative.\\
(ii) $d(u(x))=r(u(x))=x$, $\forall x\in \GR^{(0)}$, $u(r(g))g=g$ and 
$gu(d(g))=g$,  $\forall g\in \GR^{(1)}$ and $u:\Gr0 \to \Gr1$
is one-to-one.\\
(iii) $r(g^{-1})=d(g)$, $d(g^{-1})=r(g)$, $gg^{-1}=u(r(g))$ and 
$g^{-1}g=u(d(g))$.
 
\begin{definition}\label{Almost.differentiable}
 An almost differentiable groupoid 
$\GR=(\Gr0,\Gr1,d,r,\mu,u,\iota)$ is a groupoid 
such that $\Gr 0$ and $\Gr1$ are manifolds with corners, 
the structural maps $d,r,\mu,u,\iota$ are differentiable, and 
the domain map $d$ is a submersion. 
\end{definition}

We observe that $\iota$ is a diffeomorphism and hence 
$d$ is a submersion if and only if $r=d\circ \iota$ is a
submersion. Also, it follows from the definition that
each fiber $\GR_x=d^{-1}(x) \subset \Gr1$ 
is a smooth manifold whose  dimension $n$ is
constant on each connected component of $\Gr 0$. The \'etale groupoids
considered in \cite{Brylinski-Nistor1} are extreme examples
of differentiable groupoids (corresponding to
$\dim \GR_x =0$). If $\Gr0$ is smooth (\ie if it 
has no corners) then $\Gr1$ is also smooth and $\GR$ becomes
what is known as a differentiable, or Lie groupoid.\footnote{Earlier
terminology, such as in \cite{Mackenzie1}, used the 
term Lie groupoid only for differentiable groupoids in which every
pair of objects is connected by a morphism.}

We now introduce a few important geometric objects associated
to an almost differentiable groupoid.

The vertical tangent bundle (along the fibers
 of $d$)
of an almost differentiable groupoid $\GR$ is
\begin{equation}T_d \GR = \ker d_* 
= \bigcup_{x\in \Gr 0} T \GR_x \subset T\Gr 1.
\end{equation}
Its restriction  $A(\GR) = T_d \GR\big |_{\Gr0}$ to the set of units
is the Lie algebroid of $\GR$ \cite{Mackenzie1,Pradines2}. 
We denote by $T_d^* \GR$ the dual of $T_d \GR$ and by $A^*(\GR)$  
the dual of $A(\GR)$.  In addition to these bundles we shall also consider 
the bundle $\hden$ of $\lambda$-densities along the
fibers of $d$. If the fibers
  of $d$ have dimension $n$ then
$\hden=|\Lambda^n T_d^*\GR|^\lambda =\cup_x\Omega^\lambda(\GR_x).$
By invariance these bundles can be obtained as pull-backs
of bundles on $\Gr0$. For example $T_d \GR = r^*(A(\GR))$ and 
$\hden =r^*({\mathcal D}^\lambda)$ where ${\mathcal D}^\lambda$ denotes
$\hden\vert_{\Gr0}$.
If $E$ is a (smooth complex) vector bundle on the set of units
$\Gr0$ then the pull-back bundle $r^*(E)$ on $\GR$ will have 
right invariant connections obtained as follows. A connection $\nabla$
on $E$ lifts to a connection on $r^*(E)$. Its restriction to
any fiber $\GR_{x}$ defines a  linear connection  in the usual sense, 
which is denoted by $\nabla_{x}$. It is easy to see that
these connections are right invariant in the sense that
\begin{equation}
\label{eq.invariant.conn}
R_{g}^*\nabla_{x}=\nabla_{y}, \  \ \ \forall g\in \GR
\mbox{ such that } \alp (g)=x \mbox{ and }\bet (g)=y.
\end{equation}
The bundles considered above will thus have invariant connections.

The bundle $A(\GR)$, called the {\em Lie algebroid} of $\GR$,
plays in the theory of almost differentiable groupoids the r\^ole Lie 
algebras play in the theory of Lie groups. We recall for the
benefit of the reader the definition of a Lie algebroid \cite{Pradines2}.

\begin{definition}\label{Lie.Algebroid} A  Lie algebroid $A$  
over a manifold $M$ is a vector bundle $A$ over $M$ together with  a Lie
algebra structure on the space $\Gamma(A)$ of smooth sections of $A$,
and  a bundle map $\rho: A \rightarrow TP$, extended to a map between 
sections of these bundles,  such that

(i) $\rho([X,Y])=[\rho(X),\rho(Y)]$; and

(ii) $[X, fY] = f[X,Y] + (\rho(X) f)Y$

\noindent for any smooth sections $X$ and $Y$ of $A$ and any smooth function
$f$ on $M$.
\end{definition}

Note that we allow the base $M$ in the definition above to be a 
manifold with corners. 

If $\GR$ is an almost differentiable groupoid then $A(\GR)$ will
naturally have the structure of a Lie algebroid \cite{Mackenzie1}. Let
us recall how 
this structure is defined (the original definition easily
extends to include manifolds with corners). Clearly $A(\GR)$
is a vector bundle. The right translation by an arrow $g \in \GR$
defines a diffeomorphism $R_g:\GR_{r(g)}\ni g' \to g'g \in \GR_{d(g)}$. This
allows us to talk about right invariant differential geometric 
quantities as long as they are completely determined by their restriction
to all submanifolds $\GR_x$. This is true of functions and $d$--vertical
vector fields, and this is all that is needed in order to define the 
Lie algebroid structure on $A(\GR)$.
The sections of $A(\GR)$ are in
one-to-one correspondence with vector fields $X$ on $\GR$ that
are $d$--vertical, in the sense that $d_*(X(g))=0$, and right 
invariant. The condition $d_*(X(g))=0$ means that $X$ is tangent to
the submanifolds $\GR_x$, the fibers of $d$. The Lie bracket $[X,Y]$ of 
two $d$--vertical right--invariant vector fields $X$ and $Y$ will
also be $d$--vertical and right--invariant, and hence the Lie bracket 
induces a Lie algebra structure on the sections of $A(\GR)$. To define
the action of the sections of $A(\GR)$ on functions on $\Gr0$, observe that the
right invariance property makes sense also for functions  on $\GR$ and that
$\CI(\Gr0)$ may be identified with the subspace of right--invariant functions on
$\GR$. If $X$ is a right--invariant vector field  on $\GR$ and
$f$ is a right--invariant function on $\GR$ then $X(f)$ will still  be a 
right invariant function. This identifies the action of $\Gamma(A(\GR))$
on functions on $\Gr0$.

Not every Lie algebroid is the
Lie algebroid of a Lie groupoid (see \cite{al-mo:suites} for an
example).  However, every 
Lie
algebroid {\em is} associated to a {\em local} Lie groupoid 
\cite{Pradines1}. The definition
of a local Lie (or more generally, almost differentiable) groupoid 
\cite{VanEst1} is 
obtained by relaxing the condition that the multiplication $\mu$
be everywhere defined on $\Gr2$ (see Equation \eqref{Mu.everywhere}), 
and replacing it by the condition that 
$\mu$ be defined in a neighborhood $\mathcal U$ of the set of units.

\begin{definition}[van Est] 
An almost differentiable local groupoid $\mathcal L
=(\mathcal L^{(0)},\mathcal L^{(1)})$ is a pair of manifolds
with corners together with structural
morphisms $d,r: \mathcal L^{(1)} \to \mathcal L^{(0)}$, 
$\iota:\mathcal L^{(1)} \to
\mathcal L^{(1)}$, $u:\mathcal L^{(0)} \to \mathcal L^{(1)}$ and $\mu:
\mathcal U \to \mathcal L^{(1)}$, where $\mathcal U$ is a neighborhood of 
$(u \times u)(\mathcal L^{(0)})=\{(u(x),u(x))\}$
in $\mathcal L^{(2)}=\{(g,h), d(g)=r(h)\} 
\subset \mathcal L^{(1)} \times \mathcal L^{(1)}$.
The structural morphisms are required to be differentiable maps such that 
$d$ is a submersion, $u$ is an
embedding, and to satisfy the following properties:

(i) The products $u(d(g))g$, $gu(r(g))$, $g g^{-1}$ and
$g^{-1}g$ are defined and coincide with, respectively, $g$, $g$, 
$u(r(g))$ and $u(d(g))$; where we denoted $g^{-1}=\iota(g)$ as usual.

(ii) If $gh$ is defined, then $h^{-1}g^{-1}$ is defined and equal
to $(gh)^{-1}$.

(iii) (Local associativity) If $gg'$, $g'g''$ and $(gg')g''$ are defined
then $g(g'g'')$ is also defined and equal to $(gg')g''$.
\end{definition}

The set $\mathcal U$ is the set 
of arrows for which the product $gh=\mu(g,h)$ is defined.

We see that the only difference between a groupoid and a local groupoid
${\mathcal L}$ is the fact that the condition $d(g)=r(h)$ is  necessary for
the product $gh=\mu(g,h)$ to be defined, but not sufficient in general. 
The product is defined as soon as the arrows $g$ and $h$ 
are ``small enough''.
A consequence of this definition is that the right multiplication by
an arrow $g \in \mathcal L^{(1)}$ defines only a diffeomorphism 
\begin{equation}\label{local.isomorphism.for.sets}
\mathcal U_{g^{-1}} \ni g' \to g'g \in \mathcal U_{g}
\end{equation}
of an open (and possibly empty) subset $\mathcal U_{g^{-1}}$ of 
${\mathcal L}_y$, $y=r(g)$ to an  open subset 
$\mathcal U_{g} \subset {\mathcal L}_x$, $x=d(g)$. This will not
affect the considerations above, however, so  we can associate a
Lie algebroid $A(\LGR)$ to any almost differentiable local
groupoid $\LGR$. 

In the following, when considering groupoids, we shall sometimes 
refer to them as {\em global} groupoids, in order to stress the difference 
between groupoids and  local groupoids.

\section{Main definition}

Consider a complex vector bundle $E$ on the space of units $\Gr0$
of an almost differentiable groupoid $\GR$.
Denote by $r^*(E)$ its pull-back to $\Gr1$. Right translations on $\GR$
define linear isomorphisms 
\begin{gather}\label{isomorphism}
U_g:\CI(\GR_{d(g)},r^*(E)) \to \CI(\GR_{r(g)},r^*(E))  \\ \nonumber
(U_gf)(g')=f(g'g) \in (r^*E)_{g'}
\end{gather}
which  makes sense because $(r^*E)_{g'}=(r^*E)_{g'g}=E_{r(g')}$. 

If $\GR$ is merely a {\em local} groupoid then \eqref{isomorphism}
is replaced by the isomorphisms
\begin{equation}\label{Local.isomorphism}
U_g:\CI({\mathcal U}_{g},r^*(E)) \to \CI({\mathcal U}_{g^{-1}},r^*(E)) 
\end{equation}
defined for the  open subsets ${\mathcal U}_{g} \subset \GR_{d(g)}$
and ${\mathcal U}_{g^{-1}} \subset \GR_{r(g)}$ defined in 
\eqref{local.isomorphism.for.sets}.

Let $B\subset \RR^n$ be an open subset.
Define the space ${\mathcal S}^m(B \times \RR^n)$ of 
symbols on the bundle $B \times \RR^n \to B$  as in \cite{Hormander3} 
to be the set of smooth functions $a:B \times \RR^n \to \CC$ 
such that
\begin{equation}
\label{eq.symbol.estimates}
|\partial_y^\alpha \partial_\xi^\beta a(y,\xi)| \leq C_{K,\alpha,\beta}
(1+|\xi|)^{m-|\beta|}
\end{equation}
for any compact set $K\subset B$ and any multiindices $\alpha$ and $\beta$.
An element of one of our spaces $S^m$ should properly be said to have
``order less than or equal to $m$''; however, by abuse of language we
will say that it has ``order $m$''.

A symbol $a \in {\mathcal S}^m(B \times \RR^n)$ is called {\em classical} if
it has an asymptotic expansion as an infinite sum of homogeneous symbols
$a \sim \sum_{k=0}^\infty a_{m-k}$, $a_l$ homogeneous 
of degree $l$: $a_l(y,t\xi)=t^la_l(y,\xi)$ if $\|\xi\|\geq 1$
and $t \geq 1$. (``Asymptotic 
expansion''  is used here in the sense 
that $a -\sum_{k=0}^{N-1} a_{m-k}$ belongs to 
${\mathcal S}^{m-N}(B \times \RR^n)$.) The space of classical symbols
will be denoted by ${\mathcal S}^m_{\cl}(B \times \RR^n)$. We shall be
working exclusively with classical symbols in this paper.

This definition immediately extends to give spaces ${\mathcal
S}_{\cl}^m(E;F)$ of symbols on $E$ with values in $F$, where $\pi:E
\to B$ and $F \to B$ are smooth euclidian vector bundles.  These
spaces, which are independent of the metrics used in their definition,
are sometimes denoted ${\mathcal S}_{\cl}^m(E;\pi^*(F))$.  Taking $E=B
\times \RR^n$ and $F=\CC$ one recovers ${\mathcal S}_{\cl}^m(B \times
\RR^n)={\mathcal S}_{\cl}^m(B \times \RR^n;\CC)$.

A pseudodifferential operator 
$P$ on $B$ is a linear map $P:\CIc(B) \to \CI(B)$ 
that is locally of the form $P=a(y,D_y)$ plus a regularizing operator,
where for any complex valued symbol $a$ on $T^*W=W \times \RR^n$, 
$W$ an open subset of $\RR^n$, one defines $a(y,D_y):\CIc(W) \to \CI(W)$ by
\begin{equation}
a(y,D_y)u(y)=(2\pi)^{-n}\int_{\RR^n}e^{i y\cdot \xi}
a(y,\xi)\hat{u}(\xi)d\xi\,.
\end{equation} 
Recall that an operator $T:\CIc(U) \to \CI(V)$ is called
{\em regularizing} if and only if it
has a smooth distribution (or Schwartz) kernel. This happens if and only if
$T$ is pseudodifferential of order $-\infty$.

The class of $a$ in 
${\mathcal S}_{\cl}^m(T^*W)/{\mathcal S}_{\cl}^{m-1}(T^*W)$
does not depend on any choices; the collection of all these classes,
for all coordinate neighborhoods $W$, patches together to 
define a class $\sigma_m(P)\in {\mathcal S}_{\cl}^m(T^*W)
/{\mathcal S}_{\cl}^{m-1}(T^*W)$ which is called {\em the principal
symbol} of $P$. If the operator $P$ acts on sections of
a vector bundle $E$, then the principal symbol 
$\sigma_m(P)$ will belong to 
${\mathcal S}_{\cl}^m(T^*B;\End(E))/{\mathcal S}_{\cl}^{m-1}(T^*B;\End(E))$.
See \cite{Hormander3} for more details on all these constructions.

We shall sometimes refer to pseudodifferential operators acting 
on a smooth manifold as {\em ordinary} pseudodifferential operators,
in order to distinguish them from pseudodifferential operators
on groupoids, a class of operators  which we now
define (and which are really {\em families} of ordinary pseudodifferential 
operators).
 
Throughout this paper,  we shall denote  by $(P_x,x \in \Gr0)$ a family 
of order $m$ pseudodifferential operators $P_x$,
acting on the spaces $\CIc(\GR_x,r^*(E))$ for some vector bundle $E$
over $\Gr0$. Operators between sections of two {\em different} vector
bundles $E_1$ and $E_2$ are obtained by considering $E=E_1 \oplus E_2$.

\begin{definition}\label{Def.Differentiable}
A family $(P_x,x \in \Gr0)$ as above is called differentiable
if for any open set $V\subset \GR$, diffeomorphic through a 
fiber preserving diffeomorphism to $d(V) \times W$, 
for some open subset $W \subset \RR^n$, and any $\phi \in \CIc(V)$, 
we can find $a \in {\mathcal S}_{\cl}^m(d(V) \times T^*V;\End(E))$
such that $\phi P_x \phi$ corresponds to $a(x,y,D_y)$ 
under the diffeomorphism  $\GR_x\cap V \simeq W$, for each $x \in d(V)$.
\end{definition}

A fiber preserving diffeomorphism is a diffeomorphism 
$\psi:d(V) \times W \to V$ satisfying $d(\psi(x,w))=x$.
Thus we require that the operators $P_x$ be
given  in local coordinates by symbols $a_x$ that depend smoothly on
all variables, in particular 
on $x\in \Gr0$.

\begin{definition}\label{Main.definition}
An order $m$ invariant pseudodifferential operator $P$ 
on  an almost differentiable groupoid 
$\GR$, acting on sections of the vector bundle $E$,
is a differentiable family $(P_x,x \in \Gr0)$
of order $m$ classical pseudodifferential operators $P_x$
acting on $\CIc(\GR_x,r^*(E))$ and satisfying 
\begin{equation}
P_{r(g)}U_g = U_g P_{d(g)}\;(invariance)
\end{equation}
for any $g \in \Gr1$, where $U_g$ is as in \eqref{isomorphism}. 
\end{definition}

Replacing the coefficient bundle $E$ by $E \otimes \VD^\lambda$  and using the
isomorphism $\Omega_d^{\lambda}\simeq r^*(\VD^{\lambda})$, we obtain
operators acting on  
sections of density bundles.  Note that $P$ can generally {\em
not} be considered as a single pseudodifferential operator on $\Gr1$.
This is because a family of pseudodifferential operators on a smooth manifold
$M$, parametrized by a smooth manifold $B$, is not a pseudodifferential 
operator on the product $M \times B$, although it acts naturally on 
$\CIc(M \times B)$. (See \cite{Atiyah-Singer4} or \cite{Hormander3}, page 94.) 

Recall \cite{Hormander1} that distributions on a manifold $Y$ with
coefficients in the bundle $E_0$ are continuous linear
maps $\CIc(Y,E_0'\otimes \Omega) \to \CC$, where $E_0'$ is the dual
bundle to $E_0$ and  $\Omega=\Omega(Y)$
is the space of $1$-densities on $Y$. The collection of all
distributions on $Y$ with coefficients in the (finite dimensional
complex vector) bundle $E_0$ is denoted
${\mathcal C}^{-\infty}(Y;E_0)$.

If $\FAM$ is a family of pseudodifferential operators
acting on $\GR_x$ denote by $k_{x}$ the distribution kernel of $P_x$
\begin{equation}
\label{kernel}
k_{x} \in {\mathcal C}^{-\infty}(\GR_x \times \GR_x;
r_1^*(E) \otimes r_2^*(E)' \otimes \Omega_2).
\end{equation}
Here $\Omega_2$ is the pull-back of the bundle of vertical densities $\Omega_d$ 
on $\GR_x$ to $\GR_x \times \GR_x$ via the second projection. These
distribution kernels are obtained using Schwartz' kernel theorem.
We define the support of the operator $P$ to be
\begin{equation}
\label{support}\supp (P) = \overline{\cup_x \supp(k_{x})}.
\end{equation}
The support of $P$ is contained in the closed subset
$\{(g,g'), d(g)=d(g')\}$ of the product $\Gr1 \times \Gr1$. In particular
$(id \times \iota) ( \supp (P) ) \subset \Gr2$. If all operators
$P_x$ are of order $-\infty$, then each kernel $k_x$ is a smooth
section. Actually we have more

\begin{lemma}\label{Lemma.Smoothness}
The collection of all distribution kernels $k_x$ of a differentiable 
family $\FAM$ of order $-\infty$ operators 
defines a smooth section $k$ of $r_1^*(E) \otimes r_2^*(E)' \otimes \Omega_2$
on $\{(g,g'), d(g)=d(g')\}$. 
\end{lemma}

\begin{proof}Indeed if $\psi:d(V) \times W \to V$
is a fiber preserving diffeomorphism as in Definition \ref{Def.Differentiable},
then it follows from the definition that $k$ is smooth
on $d(V) \times W \times W \subset \{(g,g'), d(g)=d(g')\}$. Since
in this way we obtain an atlas of $\{(g,g'), d(g)=d(g')\}$, we obtain that
$k$ is smooth as claimed.
\end{proof}

\begin{definition}
The family $P=(P_x,x \in \Gr 0)$ is properly supported
if $p_i^{-1}(K) \cap \supp(P)$ is
a compact set for any compact subset $K \subset \GR$, where
$p_1,p_2 :\GR \times \GR \to \GR$ are the two projections.
The family  $P$ is called compactly supported if its support
$\supp(P)$ is compact; and, finally, $P$ is called uniformly supported if
its reduced support $\supp_\mu(P)=\mu_1(\supp(P))$ is 
a compact subset of $\Gr 1$, where $\mu_1(g',g)=g'g^{-1}$. 
\end{definition}

It  immediately follows from the definition that a uniformly supported 
operator is also properly supported, and that a 
compactly supported operator is uniformly supported. If the family
$\FAM$ is properly supported then each $P_x$ is properly supported, but
the converse is not true.

Recall that the composition of two ordinary 
pseudodifferential operators is defined
if one of them is properly supported. 
It follows that we can define the composition $PQ$ 
of two properly supported families of operators 
$P=(P_x,x \in \Gr0)$ and $Q=(Q_x,x \in \Gr0)$ on $\Gr 1$
by pointwise composition $PQ=(P_xQ_x,x \in \Gr0)$. The action on sections
of $r^*(E)$ is also defined pointwise as follows.
For any smooth section $f\in  \CI(\GR,r^*(E))$ denote by $f_x$ the
restriction $f\vert_{\GR_x}$. 
If each $f_x$ has compact support and $\FAM$ is a family of
ordinary pseudodifferential operators, then we define $Pf$ by 
$(Pf)_x=P_{x}(f_x).$ 

\begin{lemma}\label{Lemma.Differentiable}
(i) If $f \in \CIc(\GR,r^*(E))$ and $P=(P_x,x \in \Gr0)$ 
is a differentiable family of ordinary pseudodifferential operators 
then $Pf \in \CI(\GR,r^*(E))$. If $P$ is also properly supported then
$Pf \in \CIc(\GR,r^*(E))$.

(ii) The composition $PQ=(P_xQ_x,x \in \Gr0)$
of two properly supported differentiable families of operators 
$P=(P_x,x \in \Gr0)$ and $Q=(Q_x,x \in \Gr0)$
is a properly supported differentiable family. 
\end{lemma}

\begin{proof} If $P$ consists of regularizing operators 
then 
$$
Pf(g)=\int_{\GR_x}k_x(g,h)f(h)\, , \text{ where } x=d(g).
$$
Lemma \ref{Lemma.Smoothness} implies that the formula above for 
$Pf$ involves only the integration of smooth 
(uniformly in $g$) compactly supported
sections, and hence we can exchange integration and derivation to obtain 
the smoothness of $Pf$. This proves (i) in case
$P$ consists of regularizing operators. 
The proof of (ii) if both $P$ and $Q$ consist of regularizing
operators follows the same reasoning.

We prove now (i) for $P$ arbitrary.
Fix $g \in \GR_x$ and $V$ a neighborhood of
$g$ fiber preserving diffeomorphic to $d(V) \times W$ for some
open convex subset $W$ in $\RR^n$, $0 \in W$, such that $(x,0)$ maps to $g$. 
Replacing $P_x$ by $P_x-R_x$ for a smooth
regularizing family $R_x$ we can assume that the distribution kernels
$k_x$ of $P_x$ satisfy 
$$p_1^{-1}(d(V)\times W/4) \cap \overline{\cup \supp(k_x)}
\subset (d(V)\times W/4) \times (d(V)\times W/2).$$
The smoothness of $Pf$, respectively of $PQ$ if $Q$ consists of regularizing
operators, reduces in this way to a computation in local
coordinates. This completes the proof of (i) in general, 
and of (ii) if $Q$ is regularizing.

For arbitrary $Q$ we can replace $Q$, in view of what has
already been proved, with $Q-R$, where $R$ is a regularizing
family. In this way we may assume that
$$p_1^{-1}(d(V)\times W/2) \cap \overline{\cup \supp(k_x')}
\subset (d(V)\times W/2) \times (d(V)\times 3W/4)$$
where $k_x'$ are the distribution kernels of $Q_x$.  The support
estimates above for  $P$ and $Q$ show that the $P_yQ_y$ for $y \in d(V)$ 
are the compositions of smooth families of pseudodifferential 
operators acting on $W\subset\RR^n$. The result is then known.
\end{proof}

The smaller class of uniformly supported operators is also closed under
composition.

\begin{lemma}\label{Lemma.Support}
The composition $PQ=(P_xQ_x,x \in \Gr0)$
of two uniformly supported families of operators 
$P=(P_x,x \in \Gr0)$ and $Q=(Q_x,x \in \Gr0)$
is uniformly supported. 
\end{lemma}

\begin{proof} The reduced support $\supp_\mu(PQ)$ (see \eqref{support})
of the composition $PQ$ satisfies
$$
\supp_\mu(PQ) \subset \mu\big( \supp_\mu(P)\times \supp_\mu(Q)\big) 
$$
where $\mu$ is the composition of arrows. Since 
$\supp_\mu(P)$ and $\supp_\mu(Q)$ are compact, the  equation above 
completes the proof of  the lemma. 
\end{proof}

Let $\GR$ be an almost differentiable groupoid.
The space of order $m$, invariant, {\em uniformly} supported pseudodifferential
operators on $\GR$, acting on  sections of the vector bundle $E$
will  be  denoted by $\PS m$. We denote 
$\PS{\infty}=\cup_{m \in \ZZ}\PS{m}$
and $\PS{-\infty}=\cap_{m \in \ZZ}\PS{m}$. Thus an operator
$P \in \PS{m}$ is actually a differentiable family $\FAM$ of ordinary
pseudodifferential operators.

\begin{theorem}\label{Theorem.Algebra} 
The set $\PS{\infty}$ of uniformly supported
invariant pseudodifferential operators on an almost differentiable groupoid
$\GR$ is a filtered  algebra, \ie $$\PS{m}\PS{m'}\subset \PS{m+m'}.$$
In particular $\PS{-\infty}$ is a two-sided ideal.
\end{theorem}

\begin{proof} Let $P=(P_x,x \in \Gr0)$ and $Q=(Q_x,x \in \Gr0)$
be two invariant uniformly supported pseudodifferential operators on
$\GR$, of order $m$ and $m'$ respectively.  Their composition
$PQ=(P_xQ_x)$, is a uniformly supported operator of order $m+m'$, in
view of  Lemma \ref{Lemma.Support}.  It is also a differentiable
family due to Lemma \ref{Lemma.Differentiable}.  We now check the
invariance condition. Let $g$ be an arbitrary arrow and $U_g :
\CIc(\GR_{x},r^*(E)) \to \CIc(\GR_{y},r^*(E))$, $x=d(g)$ and $y=r(g)$,
be as in the definition above.  Then
\begin{equation*}
(PQ)_{y} U_g = P_{y} Q_{y} U_g =
P_{y}U_g Q_{x} = U_g P_{x} Q_{x} = U_g (PQ)_{x}.
\end{equation*}
This proves the theorem. 
\end{proof}

Properly supported invariant differentiable families
of pseudodifferential operators also form
a filtered algebra, denoted $\Psi_{\prop}^{\infty}(\GR;E)$.  
While it is clear that in order for our class of pseudodifferential 
operators to form an algebra we need some condition on the support 
of their distribution kernels,  exactly what support condition to 
impose is a matter of choice.
We prefer the uniform support condition because it
leads to a better control at infinity of the family of operators 
$\FAM$ and allows us to identify the regularizing ideal (\ie
the ideal of order $-\infty$ operators) with the groupoid convolution 
algebra of $\GR$. The choice of uniform support
will also ensure that $\PS{m}$ behaves functorially with respect to open
embeddings. The compact support condition enjoys the same properties
but is usually too restrictive.
The issue of support will be discussed again in examples.

The definition of the principal symbol extends easily to  $\PS{m}$.
Denote by  $\pi:A^*(\GR)\rightarrow M, \ (M=\Gr0 ) $ the projection. 
If $\FAM \in \PS{m}$ is an order $m$ 
pseudodifferential differential operator on $\GR$, then the principal
symbol $\sigma_m(P)$ of $P$ will be represented by sections of the bundle
$\End(\pi^*E)$ and will be defined to satisfy 
\begin{equation}
\label{princ.symb}
\sigma_m(P)(\xi)=\sigma_m(P_x)(\xi) \in \End(E_x) 
\; \text{ if } \xi \in A^*_x(\GR)=T^*_x\GR_x
\end{equation}
(the equation above is $\mod {\mathcal S}_{\cl}^{m-1}(A^*_x(\GR);\End(E))$). 
This equation will obviously uniquely determine a linear map 
\begin{equation*}
\sigma_m:\tPS{m} \to {\mathcal S}_{\cl}^m(A^*(\GR);\End(E))/
{\mathcal S}_{\cl}^{m-1}(A^*(\GR);\End(E)).
\end{equation*}
provided we can show that for any $\FAM$  there exists 
a symbol $a \in {\mathcal S}_{\cl}^m(A^*(\GR);\End(E))$ whose restriction to 
$A^*_x(\GR)$ is a representative of the principal symbol of $P_x$
in that fiber for each $x$. 
We thus need to choose for each $P_x$ a representative $a_x \in 
{\mathcal S}_{\cl}^{m}(A^*_x(\GR);\End(E))$ of $\sigma_m(P_x)$
such that the family $a_x$ is
smooth and invariant. Assume first that  $E$ 
is the trivial line bundle and  proceed as in
\cite{Hormander3} Section 18.1, especially Equation (18.1.27) and
below.

Choose a connection $\nabla$  on the vector bundle
$A(\GR) \to \Gr0$ and consider the pull-back vector bundle $r^{*}(A)\to \GR$
of $A(\GR) \to \Gr0$ endowed with the pull-back connection 
$\widetilde{\nabla}=r^*\nabla$.  Its  restriction
on any fiber $\GR_{x}$ defines a  linear connection  in the usual sense, 
which is denoted by $\nabla_{x}$. These connections are right invariant 
in the sense that
\begin{equation}
\label{eq.invariant.connection}
R_{g}^*\nabla_{x}=\nabla_{y}, \  \ \ \forall g\in \GR
\mbox{ such that } \alp (g)=x \mbox{ and }\bet (g)=y.
\end{equation}

Using such an invariant connection, we may define the 
exponential map of a Lie algebroid, which generalizes
the usual exponential map of  a manifold with a connection
and the exponential map of a Lie algebra as follows.
For any $x\in \Gr0$, define a   map $\exp_{x} : A_{x} \to \GR$ 
as the composition of the maps:
$$A_{x} \stackrel{i}{\longrightarrow}T_x\GR_{x}
\stackrel{\tilde{\exp}_{x}}{-\!\!\!-\!\!\!-\!\!\!\longrightarrow}\GR$$
where $i$ is the natural inclusion and
$\tilde{\exp}_{x}=\exp_{\nabla_x}$ 
is the usual exponential map  at $x\in \GR_{x}$
on the manifold  $\GR_{x}$. By varying the point $x$, we obtain
a map $\exp_{\nabla}$ defined in a neighborhood of the
zero section,  called the {\em exponential map} of the
Lie algebroid\footnote{See \cite{Landsman} for
an alternative definition of the exponential map.  One
should not confuse  this map   with the exponential
map from $\gm (A) $ to the bisections of the groupoid
as defined in \cite{KumperaS}.}.
Clearly, $\exp_\nabla$ is a local diffeomorphism 
\begin{equation}\label{eq.exp.diffeomorphism}
A(\GR) \supset V_0 \ni v \longrightarrow \exp_\nabla(v)=y \in V \subset \GR
\end{equation} 
mapping an open neighborhood $V_0$ of the zero section in $A(\GR)$ 
diffeomorphically to a neighborhood $V$ of $\Gr0$ in $\GR$, and sending 
the zero section onto the set of units.
Choose a cut-off function $\phi \in \CI(\GR)$
with support in $V$ and equal to $1$ in a smaller
neighborhood of $\Gr0$ in $\GR$.
If $y \in V$, $x=d(y)$ and $\xi \in A^*_x(\GR)$ 
let $v \in V_0$ be the unique vector $v \in A_x(\GR)$
such that $y=\exp_\nabla(v)$ and denote $e_\xi(y)=\phi(y)e^{iv\cdot \xi}$ which
extends then to all $y \in \GR$  due to the cut-off
function $\phi$. Define the $(\nabla,\phi)$--complete symbol 
$\sigma_{\nabla,\phi}(P) $ by 
\begin{equation}
\label{eq:symbol}
\sigma_{\nabla,\phi}(P)(\xi)=
(P_x e_\xi)(x), \ \ \ \ \forall \xi \in  T^*_{x}\GR_x =A^*_x(\GR).
\end{equation}

\begin{lemma}\label{lemma.symb} If $\FAM$ is an operator in $\tPS{m}$ then the
function $\sigma_{\nabla,\phi}(P)$ defined above is differentiable
and defines a symbol in ${\mathcal S}_{\cl}^m(A^*(\GR))$. Moreover
if $(\nabla_1,\phi_1)$ is another pair consisting of an invariant
connection $\nabla_1$ and a cut-off function $\phi_1$ then
$\sigma_{\nabla,\phi}(P)-\sigma_{\nabla,\phi_1}(P)$ is
in ${\mathcal S}_{\cl}^{-\infty}(A^*(\GR))$ and 
$\sigma_{\nabla,\phi}(P)-\sigma_{\nabla_1,\phi_1}(P)$ is
in ${\mathcal S}_{\cl}^{m-1}(A^*(\GR))$.
\end{lemma}

\begin{proof} For each $\xi \in A_x^*$ the function $e_\xi$ is smooth
with compact support on $\GR_x$ so $P_xe_\xi$ is defined.  Equation 
(18.1.27) of \cite{Hormander3} shows that 
$a(\xi)=\sigma_{\nabla,\phi}(P)(\xi)$ is the restriction of the 
complete symbol of $P_x\phi$ to $T^*_x\GR_x$ if the complete symbol is
defined in the normal coordinate system at $x \in \GR_x$ 
(given by the exponential map). The normal coordinate system
defines, using a local trivialization of $A(\GR)$, a fiber preserving diffeomorphism
$\psi:d(V) \times W \to V$  for some open subset $W$ of
$\RR^n$ (\ie satisfying $d(\psi(x,w))=x$). From the definition of the
smoothness of the family $P_x$ (Definition \ref{Def.Differentiable})
it follows that the
complete symbol of $P\phi$ is in $\mathcal{S}_{\cl}^m(d(V) \times T^*W)$
if the support of $\phi$ is chosen to be in $V$. This proves that
$\sigma_{\nabla,\phi}(P)$ is in ${\mathcal S}_{\cl}^m(A^*(\GR))$.

The rest follows in exactly the same way.
\end{proof}

The lemma above justifies the following definition of the principal
symbol as the class of $\sigma_{\nabla,\phi}(P)$ modulo terms of lower order
(for the trivial line bundle $E=\CC$).
This definition will be, in view of the same lemma, 
independent on the choice of $\nabla$ or $\phi$
and will satisfy Equation \eqref{princ.symb}. If $E$ is not 
trivial one can still define a complete symbol 
$\sigma_{\nabla,\nabla',\phi}(P)$, depending also on a second 
connection  $\nabla'$
on the bundle $E$, which is used to trivialize $r^*(E)$
on $V\subset \GR$ (assuming also that $V_0$ is convex). 
Alternatively, we can use Proposition \ref{prop.order.zero.red} below.

\begin{proposition} Let $\nabla$ and $\phi$ be as above.
The choice of a connection $\nabla'$
on $E$ defines a complete symbol map
$\sigma_{\nabla,\nabla',\phi}:\PS{m}\to {\mathcal S}_{\cl}^{m}(A^*(\GR)).$
The principal symbol $\sigma_m:\PS{m} \to 
{\mathcal S}_{\cl}^{m}(A^*(\GR))/{\mathcal S}_{\cl}^{m-1}(A^*(\GR))$, defined by
\begin{equation}\label{Eq.Def.PS}
\sigma_m(P) = \sigma_{\nabla,\nabla' , \phi}(P) + {\mathcal S}_{\cl}^{m-1}(A^*(\GR))
\end{equation}
does not depend on the choice of the connections $\nabla$, 
$\nabla'$ or the cut-off function $\phi$.
\end{proposition}

\begin{proof} The $(\nabla,\nabla',\phi)$--complete symbol 
$\sigma_{\nabla,\nabla',\phi}(P)$ is defined as follows. Let $w$ be a vector
in $E_x$. Using the connection $\nabla'$ we can define a section
$\tilde w$ of $r^*(E)$ on $\GR_x \cap V$ by parallel transport
along the geodesics of $\nabla$ starting at $x$, and which coincides with $w$
at $x$. Then denote $e_{\xi,w}=e_{\xi}\tilde w$
and let  
\begin{equation}
\label{eq:symbol.bundle}
\sigma_{\nabla,\nabla',\phi}(P)(\xi)w=
(P_x e_{\xi,w})(x) \in E_x, \ \ \ \ \forall \xi \in  T^*_{x}\GR_x =A^*_x(\GR).
\end{equation}
The rest of the proof proceeds along the lines of the proof of
Lemma \ref{lemma.symb}.
\end{proof}

Note that the principal symbol of $P$ determines the principal symbols of
the individual operators $P_x$ by the invariance with respect to right
translations. Precisely, we have $\sigma_m(P_x)=r^*(\sigma(P))\vert_{T^*\GR_x}$.

The following result extends some very well known properties of the calculus
of pseudodifferential operators on smooth manifolds. We shall prove
the surjectivity of the principal symbol in section \ref{D.k}.

\begin{proposition}
(i) The principal symbol map 
$$\sigma_m:\PS{m} \to {\mathcal S}_{\cl}^m(\GR;\End(E))/
{\mathcal S}_{\cl}^{m-1}(\GR;\End(E))
$$ 
has kernel $\PS{m-1}$ and satisfies Equation \eqref{princ.symb}.

(ii) The composition $PQ$ of two operators
$P,Q \in\PS{\infty}$, of orders $m$ and, respectively, $m'$,
satisfies $\sigma_{m+m'}(PQ)=\sigma_m(P)\sigma_{m'}(Q)$.
\end{proposition}

\begin{proof} (i) The operator $\FAM \in \PS{m}$ is in the kernel of
$\sigma_m$ if and only if all symbols $\sigma_m(P_x)$ vanish.
This implies $P_x \in \PS{m-1}$ for all $x$ and hence $\FAM \in \PS{m-1}$.
As already observed for $E$ a trivial line bundle,
the fact that Equation \eqref{princ.symb} is satisfied was
contained in the proof of Lemma \ref{lemma.symb}. The general case is
similar or can be proved using Proposition \ref{prop.order.zero.red}.

The second statement is known for pseudodifferential operators on
smooth manifolds \cite{Hormander3}; this accounts for the second
equality sign in the next equation. We obtain using Equation \eqref{princ.symb}
that
$$\sigma_{m+m'}(PQ)(v)=
\sigma_{m+m'}(P_xQ_x)(v)=\sigma_{m}(P_x)\sigma_{m'}(Q_x)(v)=
\sigma_{m}(P_x)(v) \sigma_{m'}(Q_x)(v)$$ where $v \in A^*_x(\GR)$.
\end{proof}

Although for the most of this paper we shall be concerned with groupoids, the 
definition of $\PS{m}$ easily extends to {\em local} groupoids. Indeed it
suffices to modify the invariance condition  in Definition
\ref{Main.definition}, using the notation in Equation \eqref{Local.isomorphism},
as follows. We assume that for any $g \in \Gr 1$ and any smooth
compactly supported function $\phi$ on ${\mathcal U}_{g}$ there
exists a regularizing operator $R_{g,\phi}$ such that 
\begin{equation}\label{Local.invariance}
U_g(\phi)P_{r(g)}U_g f - U_g (\phi P_{d(g)}f) =R_{g,\phi}f
\end{equation}
for any function $f \in \CIc({\mathcal U}_{g})$.
We thus replace the strict invariance of the original definition
by `invariance up to regularizing operators'. 

We denote by $\Psi_{\loc}^m(\GR;E)$
the set of differentiable properly supported 
families $\FAM$ of order $m$ pseudodifferential operators
satisfying the condition \eqref{Local.invariance} above. Note that if we regard an 
almost differentiable groupoid $\GR$ as a local groupoid then
$\PS{\infty} \subset \Psi_{\loc}^m(\GR;E)$. The inclusion is generally a
strict one, though, because Equation
\eqref{Local.invariance} gives no condition for order $-\infty$
operators, and so $\Psi_{\loc}^{-\infty}(\GR;E)$ consists of arbitrary smooth
families $\FAM$ of regularizing operators. This ideal is too
big to reflect the structure of $\GR$. The 
``symbolic'' part remains  however the same:
$$\Psi_{\prop}^{m}(\GR;E)/\Psi_{\prop}^{-\infty}(\GR;E) \simeq 
\Psi_{\loc}^m(\GR;E)/\Psi_{\loc}^{-\infty}(\GR;E).$$

If the sets ${\mathcal U}_g$ are all connected (in which case
the local groupoid $\GR$ is said to be $d$--connected) an easier
condition to use than \eqref{Local.invariance} is
\begin{equation}\label{Infinitesimal.invariance}
[X, P] \in \Psi_{\loc}^{-\infty}(\GR;E)
\end{equation}
for all $r$--vertical left--invariant vector fields $X$ on $\Gr 1$. With this, 
the following analog of Theorem \ref{Theorem.Algebra}, 
becomes straightforward.

\begin{theorem}\label{local.algebra} 
Assume that $\GR$ is a $d$--connected almost differentiable groupoid.
Then the space  $\Psi_{\loc}^{\infty}(\GR)$ is a filtered  algebra, 
with $\PS{-\infty}$ as residual ideal.
\end{theorem}

\begin{proof}The only thing to check is that $\ttPS{\infty}$ is closed
under composition. The composition of two differentiable, properly
supported families $P, Q \in \ttPS{\infty}$
is again differentiable and properly supported, as
it has already been proved. The infinitesimal invariance condition
$[X,PQ]=[X,P]Q + P[X,Q] \in\ttPS{-\infty}$ \eqref{Infinitesimal.invariance}
follows from the fact that
$\ttPS{-\infty}$ is an ideal of $\ttPS{\infty}$.
\end{proof}

\section{Differential operators and quantization}

In this section, we examine the differential operators in $\PS{\infty}$,
if $\GR$ is a global groupoid, or in  $\ttPS{\infty}$
if $\GR$ is a local groupoid. We also show how a simple algebraic 
construction applied to $\GR$ and to the algebras
$\ttPS{\infty}$ leads to a concrete  
construction of a deformation quantization of the 
Lie-Poisson structure on the dual of a Lie algebroid.

In this section, $\GR$ will be an almost differentiable local groupoid,
This generality is necessary in order to integrate arbitrary Lie
algebroids. Nevertheless, when
$A$ is the Lie algebroid of an almost differentiable {\em global}
groupoid $\GR$ (that is not just a {\em local} groupoid) then all
results we shall prove for $\ttPS{\infty}$ in this section
extend immediately to $\PS{\infty}$, although
we shall not mention this each time.

\begin{lemma}\label{prop.order.zero}
Let $\FAM$ be an operator in $\ttPS{\infty}$. If 
$P_x$ is a multiplication operator for all $x$, then there exists
a smooth endomorphism $\A$ of $E$ such that $P_x(g) = \A(r(g))$
for all $g \in \GR_x$. Conversely, every smooth section $\A$
of $\End(E)$ defines a multiplication operator in $\ttPS0$.
\end{lemma}

\begin{proof} By assumption $P_x(g)$ is in 
$\End(E_{r(g)})$. The invariance relation shows that $P_x(g)$ depends only
on $r(g)$. This defines the section $\A$ of $\End(E)$ such that $P_x(g) = \A(r(g))$.
To show that $\A$ is smooth, we let $\phi$ be a smooth section of $E$ 
over $\Gr0$ and let 
$\tilde{\phi}(g)=\phi(r(g))$. By assumption $P\tilde{\phi}$ is smooth
and hence $\A\phi=P\tilde{\phi}|_{\Gr0}$ is also smooth.
Since $\phi$ is arbitrary this implies the smoothness of $\A$.

Conversely, if $\A$ is a smooth endomorphism of $E$, then if we let
$P_x(g)=\A(r(x))$ we obtain a multiplication operator in $\ttPS0$.
\end{proof}

The following proposition will allow us to assume that 
$E$ is a trivial bundle, which is sometimes useful in 
applications.

\begin{proposition}\label{prop.order.zero.red} 
Let $E$ be a vector bundle on $\Gr0$ embedded
into a trivial hermitian bundle, $E \subset \CC^N$. Denote by
$e_0$ the projection onto $E$ regarded as a matrix of multiplication
operators in $M_N(\tttPS{0})$, 
the algebra of $N\times N$ matrices with values in
$\tttPS{0}$. Then  
$\ttPS{\infty} \simeq e_0M_N(\tttPS{\infty})e_0$ as filtered algebras.
\end{proposition}

\begin{proof} The multiplication operator $e_0$ defines an element
of $\ttPS0$ by the lemma above; hence it acts on all spaces 
$\CIc(\GR_x,\CC^N)$. Then 
$$\CIc(\GR_x,r^*(E))=e_0\CIc(\GR_x,\CC^N)$$
and every pseudodifferential
operator $P_x$ on $\CIc(\GR_x,r^*(E))$ extends in this way to an operator on 
$\CIc(\GR_x,\CC^N)$. This gives an inclusion
$\ttPS{\infty} \subset e_0M_N(\tttPS{\infty})e_0$. Conversely if $P_x$
is a pseudodifferential
operator on $\CIc(\GR_x,\CC^N)$ then $e_0P_xe_0$ is a 
pseudodifferential  operator on $\CIc(\GR_x,\CC^N)$.
This gives the opposite inclusion.
\end{proof}

The following proposition shows the intimate connection between
$A(\GR)$,  the Lie algebroid of $\GR$,  and $\tttPS{\infty}$
($\tPS{\infty}$ if $\GR$ is global). It is morally an equivalent 
definition of the Lie algebroid
associated to an almost differentiable local groupoid.

\begin{proposition}\label{Prop.inclusion} Let $\GR$ be an
almost differentiable local groupoid. 

(i) The algebra $\CI(\Gr0)$ is the
algebra of multiplication operators in $\tttPS{0}$.

(ii) The space of sections of the Lie algebroid $A(\GR)$ can be identified 
with the space of order $1$ differential operators in $\tttPS{1}$ without
constant term. 

(iii) The Lie algebroid structure of $A(\GR)$ is induced by the commutator
operations $[\;,\,]:\tttPS{1}\times \tttPS{1}\to \tttPS{1}$ and
$[\;,\,]:\tttPS{1} \times \tttPS{0}\to \tttPS{0}$.
\end{proposition}

\begin{proof} The first part is a particular case of Lemma
\ref{prop.order.zero}, only easier. Order $1$ differential operators 
without constant term are vector fields,  right invariant by the
definition of $\tttPS{1}$, so they can be
 identified with the sections of the
Lie algebroid $A(\GR)$ of $\GR$. This proves (ii). In order to 
check (iii) recall that, if we regard vector fields on $\GR$ as 
linear maps $\CI(\GR) \to \CI(\GR)$, then the Lie bracket coincides
with the commutator of linear maps. Moreover the commutator $[X,f]$ of
a vector field $X$ and of a multiplication map $f$ is $[X,f]=X(f)$, again
regarded as a linear map. Then (iii) follows in view of the discussion above. 
\end{proof}

The Lie algebroid $A=A(\GR)$ turns out to determine the structure of the
algebra of invariant tangential differential operators on $\GR$,
denoted $\operatorname{Diff}(\GR)$.
We shall see that the subalgebra 
$\operatorname{Diff}(\GR)\subset\tttPS{\infty}$ is a concrete 
model of the {\em  universal enveloping algebra} of the Lie algebroid $A$ 
\cite{Huebschmann,Rinehart}, a concept whose definition we now recall.

Given a Lie algebroid $A\to M$ with anchor $\rho$,  we can
make the $C^{\infty}(M)$-module direct sum $C^{\infty}(M)\oplus \gm(A)$
into a Lie algebra over $\CC$ by defining
$$[f+X, \ g+Y]=(\rho(X ) g -\rho(Y) f)+[X, Y].$$
Let $U=U(C^{\infty}(M)\oplus \gm(A))$ be its  universal enveloping algebra.  
For any $f\in C^{\infty}(M)$ and $X\in \gm (A)$, denote
by $f'$ and $X'$ their canonical image in  $U$.
Denote by $I$ the two-sided  ideal  of $U$  generated by all 
elements of the form $(fg)'-f'g'$ and $(fX)'-f'X'$.
Define
\begin{equation}
U(A)=U/I.
\end{equation}
$U(A)$ is called the universal enveloping algebra of the Lie algebroid $A$.
When $A$ is a Lie algebra, this definition reduces to the  usual universal 
enveloping algebra. We shall see, for example, that for the tangent
bundle $TM$ this is the algebra of differential operators on $M$.

The maps $f \to f'$ and $X \to X'$ considered above descend to 
linear embeddings  $i_{1}: C^{\infty}(M)\to U(A)$, and
$i_{2}: \gm (A)\to U(A)$; the first map $i_1$ is an 
algebra morphism. These maps have the following properties:
\begin{equation}
\label{eq.relations}
i_{1}(f)i_{2}(X)=i_{2}(fX), \ \ 
[i_{2}(X),i_{1}(f)]=i_{1}(\rho(X)f)), \ \
[i_{2}(X),i_{2}(Y)]=i_{2}([X,Y]).
\end{equation}
In fact, $U(A)$ is universal among triples $(B, \phi_{1}, \phi_{2})$
having these properties (see \cite{Huebschmann} for a proof of this
easy fact).

In particular, if $M$ is the space of units of an almost differentiable
groupoid $\GR$ with Lie algebroid $A(\GR)$, 
the natural morphisms $\phi_{1}:  C^{\infty}(M)
\to  \operatorname{Diff}(\GR)$ 
and $\phi_{2}: \gm (A)\to  \operatorname{Diff}(\GR)$
obtained from  Proposition \ref{Prop.inclusion} extend to 
a  unique algebra morphism $\tau :U(A) \to \operatorname{Diff}(\GR)$.
(Recall that we denoted by $\operatorname{Diff}(\GR)$
the algebra of right invariant tangential differential operators on $\GGR$.) 
Denote by $U_n(A) \subset U(A)$ the space 
generated by $\CI(M)$ and the images of
$X_1 \otimes X_2 \otimes \ldots \otimes X_k \in U=U(C^{\infty}(M)\oplus \gm(A))$,
for $k\leq n$, under the canonical projection $U \to U(A)=U/I$. Then
\begin{equation}
\label{eq.filtration}
U_{0}(A)\subset U_{1}(A) \subset \cdots \subset U_{n}(A) \subset \cdots 
\end{equation}
is a filtration of $U(A)$. The relations \eqref{eq.relations}
show that, as in the Lie algebra case, 
the graded algebra $\oplus U_n(A)/U_{n-1}(A)$ is commutative. Similarly, 
$\operatorname{Diff}(\GR)$ is naturally filtered by degree.
\begin{lemma}
\label{pro:uadiff}
The map $\tau :U(A) \to \operatorname{Diff}(\GR)$ maps
$U_n(A)$ onto the space $\operatorname{Diff}_n(\GR)$ of operators of
order $\leq n$. 
\end{lemma}

\begin{proof}
Let $D \in \operatorname{Diff}(\GR)$ be an invariant tangetial differential
operator of order $\leq n$.
By right invariance  $D$ is completely determined by the restrictions
$(Du)\vert_{\Gr0}$, $u \in \CIc(\GR)$. Since $D$ acts on the fibers of $d$ we can write
$$
(Du)\vert_{\Gr0}=\sum_{i=1}^n D_iu
$$
where $D_i$ is a superposition of derivations
$D_iu=X_1^{(i)}X_2^{(i)}\cdots X_{k_i}^{(i)}u$ defined 
using the tangential derivations $X_j^{(i)}\in \Gamma(A)$.
By definition it follows that $D$ is the sum of 
$\tau(X_1^{(i)}X_2^{(i)}\cdots X_{k_i}^{(i)})$.
\end{proof}

Denote by $\operatorname{Symm}(A)$ the {\em symmetric} tensor product
of the bundle $A$, that is 
$$\operatorname{Symm}(A)=\bigoplus_{n=0}^\infty S_n(A)$$
where $S_n(A)$ is the symmetric quotient of the 
bundle $A^{\otimes n}$, and is isomorphic
to the subspace of symmetric tensors, if $S_0(A)$
is $\CI(M)$ by convention. The space $\Gamma(\operatorname{Symm}(A))$
of smooth sections of $\operatorname{Symm}(A)$ identifies with
the space of smooth functions on $A^*$ polynomial in each fiber.
The complete symbol map $\sigma_{\nabla,\phi}(D)$ of an invariant 
differential operator $\operatorname{Diff}(\GR) \subset \tttPS{\infty}$ 
(defined in Equation \eqref{eq:symbol}) does not depend on
the cut-off function $\phi$ and will be a polynomial in $\xi$,
denoted simply by $\sigma_{\nabla}(D)$.

Using the algebra morphism $\tau :U(A) \to \operatorname{Diff}(\GR)$
obtained from the universality property of $U(A)$, we have the 
following Poincar\'e-Birkhoff-Witt type theorem for Lie algebroids.
Note that both $U(A)$ and $\gm (\operatorname{Symm}(A))$  have 
natural filtrations (see \eqref{eq.filtration}). 

\begin{theorem}[Poincar\'e-Birkhoff-Witt]
\label{Theorem.PBW}
The composite map
$$U(A) \ni D \to \sigma_{\nabla}(\tau(D)) \in \Gamma(\operatorname{Symm}(A))$$
is an isomorphism of filtered vector spaces. In particular
$\tau:U(A) \to \operatorname{Diff}(\GR)$ is an algebra isomorphism.
\end{theorem}

\begin{proof} It follows from definitions that the
map $\sigma=\sigma_{\nabla}\circ \tau$ considered in the statement
maps $U_n(A)$ to $\oplus_{k=0}^n\Gamma(S_k(A))$ and hence it
preserves the filtration. By abuse of notation we shall still denote by 
$\sigma$ the induced map $U_n(A)/U_{n-1}(A) \to \Gamma(S_n(A))$.
It is enough to prove that the map of graded spaces 
$$\sigma:\bigoplus U_n(A)/U_{n-1}(A) \to \bigoplus \Gamma(S_n(A))=
\Gamma(\operatorname{Symm}(A))$$ is
an isomorphism. By Lemma \ref{pro:uadiff} this map is onto. We now prove 
that it is one-to-one.

The inclusion of $\CI(M)$ in $U(A)$ makes
$U(A)$ a $\CI(M)$--bimodule. The filtration $U_n(A)$ of $U(A)$ consists
of $\CI(M)$--bimodules. Moreover, since the graded algebra
$\oplus U_n(A)/U_{n-1}(A)$ is commutative, 
the quotient $U_n(A)/U_{n-1}(A)$  consists of
central elements for this action (\ie the left and right $\CI(M)$--module 
structure coincide). It follows from the definition that the subspace
$\Gamma(A)^{\otimes n}$ of the universal enveloping
algebra $U=U(C^{\infty}(M)\oplus \gm(A))$
maps onto $U_n(A)/U_{n-1}(A)$. The previous discussion
shows that this map descends to a map from the tensor product 
$\Gamma(A )\otimes_{\CI(M)}\cdots \otimes_{\CI(M)} \Gamma(A )$
of $\CI(M)$-modules.
By the commutativity of the graded algebra of $U(A)$ this further
descends to a $\CI(M)$--linear surjective map  
$q:\Gamma(S_n(A)) \to U_n(A)/U_{n-1}(A)$. 

The composition $\sigma \circ q:\Gamma(\operatorname{Symm}(A))
\lon \Gamma(\operatorname{Symm}(A))$ is multiplicative since
both $q$ and $\sigma$ are multiplicative. Moreover $\sigma \circ q$ is
the identity when restricted to $\CI(M)$ (the order $0$ elements)
and $\Gamma(A)$ (the elements of order $1$). Since these form
a system of generators of the commutative algebra 
$\Gamma(\operatorname{Symm}(A))$
it follows that $\sigma \circ q$ is the identity. This completes the proof.
\end{proof}

{\bf Remark} The Poincar\'e-Birkhoff-Witt theorem was proved  in the  
algebraic context by Rinehart \cite{Rinehart} for $(L, R)$-algebras 
(an algebraic version
of Lie algebroids).  It essentially stated that the associated
graded algebra $gr U(A)=\oplus_{n} U_{n+1}(A) /U_n(A)$ is
isomorphic to symmetric algebra $S(\gm (A))=\Gamma(\operatorname{Symm}(A))$. 
The role of the connection $\nabla$ on $A\to M$ is to   establish
an {\em explicit} isomorphism $\sigma_\nabla\circ \tau$ between $U(A)$ and 
$\Gamma(\operatorname{Symm}(A))$.

We will now use the results of this and the previous section to
obtain an explicit deformation quantization of $A^*$. In order to do that
we need to establish the relation between commutators and the 
Poisson bracket in our calculus. 

For any $x\in  \Gr 0$,  $T^{*} \GR_{x} $ is a symplectic manifold,
so $T^{*}_{d}\GR \stackrel{def}{=}\cup_{x\in \Gr 0}T^{*} \GR_{x} $
is a regular Poisson manifold with the leafwise symplectic structures.
Now  the Poisson structure on $A^*$
can be considered as  being induced from that 
on $T^{*}_{d}\GR$. More precisely, let $\Phi : T^{*}_{\bet}\GGR \to A^*$ 
be the natural projection induced by the right translation, used to 
define a map $\Phi^*:\CI(A^*(\GR)) \to \CI(T^*_d\GR)$. We then have:

\begin{lemma}
\label{lem:phi}
The map $\Phi$ is a Poisson map.
\end{lemma}

Of course this lemma is really the definition of the Poisson
structure on $A^*$. The point is to show that the subspace
$\Phi^*(\CI(A^*(\GR)))$ of $\CI(T^*_d\GR)$ is closed under the 
Poisson bracket.

\begin{proof} 
It is enough to check that 
\begin{equation}
\label{eq.Poisson}
\Phi^*(\{f,g\})=\{\Phi^*(f),\Phi^*(g)\}
\end{equation} 
where $f$ and $g$ are two
smooth function on $A^*$ with {\em polynomial} restrictions
on each fiber of $A^*$, that is for  $f$ and $g$ in 
$\Gamma(\operatorname{Symm}(A))$. Since the Poisson bracket is a 
derivation in  each variable it is further enough to check this
for $f$ constant or $f$ linear in each fiber. If both $f$ and $g$ are
constant in each fiber then both sides of Equation
\eqref{eq.Poisson} vanish. If $f$ and $g$ are of degree one in each
fiber then they correspond to  sections $X$ and $Y$ of $A$,
and their Poisson bracket will identify to $[X,Y]$ (so in particular will
also be of degree one in each fiber and this justifies the name 
of Lie-Poisson structure for this Poisson structure). For this situation
the relation \eqref{eq.Poisson} follows from the identification of
$\Gamma(A)$ with $d$--vertical right invariant vector fields on $\GR$
and the fact that $T^*_d\GR$ is a Lie-Poisson manifold itself. The remaining
case is treated similarly.
\end{proof}

We shall use the following general fact on the principal symbols of
commutators.

\begin{proposition}
When $E$ is the trivial line bundle, the commutator $[P,Q]$ satisfies
\\$\sigma_{m+m'-1}([P,Q])=\{\sigma_m(P),\sigma_{m'}(Q)\}$,
where $\{\;,\,\}$ is the Poisson structure on $A^*(\GR)$.
\end{proposition}

\begin{proof} The map $\Phi^*:\CI(A^*(\GR)) \to \CI(T^*_d\GR)$ 
is a Poisson map according to Lemma \ref{lem:phi}. Since 
$\Phi^*(\sigma_m(P))=\sigma_m(P_x)$ on $T^*\GR_x$ the result follows from 
\begin{gather*}
\Phi^*(\sigma_{m+m'-1}([P,Q]))=\sigma_{m+m'-1}([P_x,Q_x])=
\{\sigma_m(P_x),\sigma_{m'}(Q_x)\}\\
=\{\Phi^*(\sigma_m(P)),\Phi^*(\sigma_{m'}(Q))\}=
\Phi^*(\{\sigma_m(P),\sigma_{m'}(Q)\}).
\end{gather*}
Since $\Phi^*$ is one-to-one this proves the last statement.
\end{proof}

We now use the results above to construct deformation quantizations.

Deform the Lie bracket structure on the Lie algebroid $A$ on
$M$ to obtain a new algebroid, the {\em adiabatic} algebroid $A_t$ 
associated to $A$,  defined over $M \times [0,\infty)$
as follows. As a bundle $A_t$ is the lift of the bundle 
$A$ to $M \times [0,\infty)$. Regard the sections $X$ of 
$A_t$ as functions $X:[0,\infty) \to \Gamma(A)$, $t \to X_t$. Then the 
algebroid structure is obtained by letting
\be
[X,Y]_t=t[X_t,Y_t]\ \ \ \mbox{ and }
\ee
\be
\rho(X)_t=t\rho(X_t)
\ee
so that $\rho(X)f$ is the function whose restriction to
$\{t\} \times M$ is $t\rho(X_t)(f_t)$ where for any 
$f \in \CI(M \times [0,\infty))$ we denote by $f_t \in \CI(M)$  
the restriction of $f$ to $\{t\} \times M \equiv M$.

Observe that $\CI([0,\infty))\subset \CI(M \times [0,\infty))$
is acted upon trivially by $\Gamma(A_t)$ and hence will define a central
subalgebra of the universal enveloping algebra $U(A_t)$ of the adiabatic 
Lie algebroid $A_t$. Denote by $t \in  \CI([0,\infty))$
the identity function.

\begin{theorem}
\label{Theorem.Quantization}
The inverse limit $\projlim U(A_t)/t^nU(A_t)$
is a deformation  quantization of $\Gamma(\operatorname{Symm}(A))$,
 the algebra of polynomial functions on $A^*$. Therefore,
it induces a $*$-product on the Lie-Poisson space $A^{*}$ in the sense
of \cite{BFFLS}.
\end{theorem}

\begin{proof} It follows from the PBW theorem for Lie algebroids
(Theorem \ref{Theorem.PBW}) that  the inverse limit 
$\projlim U(A_t)/t^nU(A_t)$ 
is isomorphic to $\Gamma(\operatorname{Symm}(A))[[t]]$ as a $\CC[[t]]$ 
module via the complete symbol map $\sigma_\nabla=\sigma_{\nabla,\phi}$
defined in Equation \eqref{eq:symbol}. Denote by $\{\;,\,\}'$
the Poisson bracket on $A^*_t$ and identify $\CI(A^*)$ with
the subset of functions on $A^*_t$ that do not depend on $t$.
Then $\{f,g\}'=t\{f,g\}$ if $f,g$ are smooth functions on
$A^*_t$ and $\{\;,\,\}$ is the Poisson bracket on $A^*$.

For any polynomial function $f$ on $A^*$ denote by $q(f)\in U(A_t)$ 
the element with complete symbol $\sigma_\nabla(q(f))(\xi,t)=f(\xi)$ 
obtained, as an application of the isomorphism in the PBW theorem 
for $A_t$. (We treat $\tau$ as the identity, which
justifies replacing $\sigma_\nabla \circ \tau$ with $\sigma_\nabla$.)
The proof will be complete if we check the following 
quantization relation
\begin{equation}
\label{eq.quantization}
q(f)q(g)-q(g)q(f)= t q(\{f,g\}) + t^2 h, \ \ h \in U(A_t).
\end{equation}
It is actually enough to do so for $f$ and $g$ among 
a set of generators of the algebra $\Gamma(\operatorname{Symm}(A))$.
Choose the set of generators to be the union of $\CI(M)$ and 
$\Gamma(A)$. Then Equation  \eqref{eq.quantization} will
obviously be satisfied for $f$ and $g$ in this generating set 
(with no $t^2$--term) in view of the definition of the Lie bracket 
on $A_t$.\end{proof}

{\bf Remark} When $M$ is a Lie group $G$ and $\nabla$ is the right invariant
trivial connection making all right invariant vector fields parallel,
this construction, restricted to right invariant
differential operators, reduces to the symmetrization correspondence
between  $U(\frakg )$ and $S(\frakg )$ studied by Berezin
\cite{Berezin} and Gutt \cite{gu:explicit}. See also Rieffel's paper
\cite{Rieffel2}. On the other hand, when the Lie algebroid $A$ is
the tangent bundle Lie algebroid $TP$, this construction
gives rise to a quantization for the canonical symplectic structure
on  cotangent bundle  $T^{*}M$.

The quantization of the Lie-Poisson structure on
$A^*$ was investigated by Landsman in terms of Jordan-Lie algebras 
\cite{Landsman}.  His quantization axioms are
closer to those of Rieffel's strict deformation quantization.
It was conjectured  in \cite{Landsman} that the quantization 
of $A^*$ is related to the groupoid $C^*$-algebra of the corresponding
groupoid $\GR $, and the transitive case was proved in \cite{Landsman1}.

\section{Examples}
As anticipated in the introduction, we recover many previously
defined classes of operators as pseudodifferential operators
on groupoids. We begin by showing that pseudodifferential
operators on a manifold, in the classical sense,
are obtained as a particular case of our construction. In this 
section we will consider only operators with coefficients in the
trivial line bundle  $E=\CC$. We include the description of the Lie algebroids
associated to each example.

Denote by $\Psi^m_{\prop}(M)$ the space of {\em properly supported} 
pseudodifferential operators on a smooth manifold $M$, and by
$\Psi_{\comp}^m(M)$ the subspace of operators with compactly supported 
Schwartz kernel, regarded as a distribution on $M \times M$. 

\begin{example}
Let  $M$ be a smooth manifold and $\GR=M \times M$ be the pair groupoid: 
$\Gr1=M \times M$, $\Gr0=M$,
$d(x,y)=y$, $r(x,y)=x$, $(x,y)(y,z)=(x,z)$. According to
the definition, a pseudodifferential operator $P\in \tPS{m}$
is a uniformly supported invariant {\em family} of 
pseudodifferential operators $P=(P_x,x \in M)$ on 
$M\times \{x\}$. The action by right translation with $g=(x,y)$ identifies
$M\times \{x\}$ with $M\times \{y\}$.
After we identify all fibers
 with $M$, the invariance condition
reads  $P_x=P_y$ for all $x,y$ in 
$M$. This shows that the family $P=(P_x)_{x \in M}$ is constant, and hence 
reduces to {\em one} operator $P_0$ on $M$. The family $P=(P_x)_{x \in M}$ 
is  uniformly supported
if and only if the distribution kernel of $P_0$ is compactly supported.
The family $P$ is properly supported if and only if $P_0$ is
properly supported. If $M$ is not compact then $P=(P_x)_{x \in M}$ 
will not be compactly supported
unless it vanishes. We obtain $\tPS{m}=\Psi_{\comp}^m(M)$.

In this case, the Lie algebroid  $A(\GR)$ is the tangent bundle 
$TM$. \end{example}

\begin{example} If $\GR$  has only one unit, \ie if $\GR=G$, a Lie group,
then $\tPS{m} \simeq \Psi^m_{\prop}(G)^G$, the algebra of properly
supported pseudodifferential operators on $G$, invariant
with respect to right translations. In this example, every invariant
properly supported operator is also uniformly supported. Again,
there are no nontrivial compactly supported operators unless $G$ is compact.

In this example $A(\GR)$ is the Lie algebra of $G$. 
\end{example}

We continue with some more elaborate examples.

\begin{example} If $\GR$ is the holonomy groupoid \footnote{The
holonomy groupoids of some foliations are non-Hausdorff manifolds.  We
believe that our constructions will extend to this case with the use
of the technique in \cite{co:survey} (page 564), where the groupoid
algebra is generated by continuous functions supported on Hausdorff
open sets.} of a foliation 
${\mathcal F}$ on a smooth manifold $M$, 
then $\tPS{\infty}$ is the algebra of 
pseudodifferential operators along the leaves of the foliation
\cite{Connes3,Connes1,Moore-Schochet1,Winkelnkemper1}. Suppose for
simplicity that the foliation is given by a (right) locally free 
action of a Lie group 
$G$ on a manifold $M$, and that the isotropy representation  of 
$G_x$, the stabilizer of $x$, on $N_x$, the normal space to the
orbit through $x$, is faithful. This is equivalent to the condition that the 
holonomy of the leaf passing through $x$ be isomorphic to the
discrete group $G_x$. Then the holonomy groupoid of this foliation 
is the transformation groupoid $\Gr1=X \times G$, $\Gr0=X$, 
$d(x,g)=xg$, $r(x,g)=x$
and $(x,g)(xg,g')=(x,gg')$. The algebra $\PS\infty$ consists of
families of pseudodifferential operators on $G$ parametrized by $X$,
invariant with respect to the diagonal action of $G$ and with support contained
in a set of the form $\{(x,kg,xk,g)\} \subset (X \times G)^2$ where $g \in G$
is arbitrary but $x\in X$ and $k\in G$ belong to compact sets that depend
on the family $P$.

The Lie algebroid  $A=A(\GR)$ is the  integrable 
subbundle of $TM$ corresponding to the foliation
${\mathcal F}$.
\end{example}

\begin{example} Let $\GR$ be the fundamental groupoid of a 
compact smooth manifold $M$ with fundamental group
$\pi_1(M)=\Gamma$. Recall that if we denote by $\tilde{M}$ a universal
covering of $M$ and let $\Gamma$ act by covering transformations, then
$\Gr0=\tilde{M}/\Gamma = M$, $\Gr1=(\tilde{M} \times
\tilde{M})/\Gamma$ and $d$ and $r$ are the two projections. Each fiber
$\GR_x$ can be identified with $\tilde{M}$, uniquely up to the action
of an element in $\Gamma$.  Let $P=(P_x,x \in M)$ be an invariant,
uniformly supported, pseudodifferential operator on $\GR$. Then each
$P_x$, $x \in M$ is a pseudodifferential operator on $\tilde{M}$. The invariance
condition applied to the elements $g$ such that $x=d(g)=r(g)$ implies that
each operator $P_x$ is invariant with respect to the action of
$\Gamma$. This means that we can identify $P_x$ with an operator on
$\tilde M$ and that the resulting operator does not depend on the
identification of $\GR_x$ with $\tilde M$.  Then the invariance
condition applied to an arbitrary arrow $g \in \Gr1$ gives that all operators
$P_x$ acting on $\tilde M$ coincide. We obtain $\tPS{m} \simeq
\Psi^m_{\prop}(\tilde{M})^\Gamma$, the algebra of properly supported $\Gamma$-invariant
pseudodifferential operators on the universal covering $\tilde{M}$ of
$M$. An alternative definition of this algebra using 
crossed products is given in \cite{Nistor4}. See also \cite{Connes3}.

The Lie algebroid is $TM$, as in the first example.
\end{example}

\begin{example} Let $\Gamma$ be a discrete group acting from the right by 
diffeomorphisms on a smooth compact manifold $M$. Define $\GR$ as follows,
$\GR^{(0)}=M$, $\GR^{(1)}=M\times M \times \Gamma$ with 
$d(x,y,\gamma)=y\gamma$,
$r(x,y,\gamma)=x$ and 
$(x, y, \gamma)(y\gamma,y',\gamma')=(x,y'\gamma^{-1},\gamma\gamma')$.
Then $\tPS{\infty}$ is the algebra  generated 
by $\Gamma$ and $\Psi^{\infty}_{\prop}(M)$ acting on 
$\CI(M) \otimes \CC[\Gamma]$, where $\Gamma$ acts diagonally, 
$\Psi^{\infty}_{\prop}(M)$ acts on the first variable, and
$\CC[\Gamma]$ denotes the set of finite sums of elements in $\Gamma$ with 
complex coefficients.
This algebra  coincides with the crossed product algebra 
$\Psi^{\infty}_{\prop}(M) \rtimes \Gamma=\{\sum_{i=0}^n P_i g_i, 
P_i \in \Psi^\infty_{\prop}(M), g_i \in \Gamma\}$. 
The regularizing algebra $\tPS{-\infty}$
is isomorphic to $\Psi^{-\infty}_{\prop}(M) \rtimes \Gamma \simeq 
\Psi^{-\infty}_{\prop}(M) \otimes \CC[\Gamma]$.
If we drop the condition that $M$ be compact we obtain 
$\tPS{\infty} \simeq \Psi_{\comp}^{\infty}(M) \rtimes \Gamma$.

In general if a discrete group $\Gamma$ acts on a groupoid $\GR_0$ then 
$$\Psi^\infty_{\prop}(\GR_0 \rtimes \Gamma) \simeq 
\Psi^\infty_{\prop}(\GR_0) \rtimes \Gamma.$$

This construction does not change the Lie algebroid.
\end{example}

In the following example we realize the algebra of
families of operators in $\Psi^m(\GR)$ parametrized by a compact
space $B$ as the algebra of pseudodifferential operators on 
the product groupoid $\GR \times B$. 
This example shows that our class of operators on groupoids is closed 
under formation of families of operators.

\begin{example} If $B$ is a compact manifold with
corners,  
define $\GR \times B$ by $(\GR \times B)^{(0)}=\Gr0 \times B$,
$(\GR \times B)^{(1)}=\Gr1 \times B$ with the structural maps preserving
the $B$-component. Then $\Psi^{m}(\GR \times B)$ contains
$\Psi^m(\GR) \otimes \CI(B)$ as a dense subset in the sense
that $\Psi^{m-1}(\GR) \otimes \CI(B) =\Psi^{m}(\GR) \otimes \CI(B)
\cap \Psi^{m-1}(\GR \times B)$ and $\Psi^{m}(\GR) \otimes \CI(B)/
\Psi^{m-1}(\GR) \otimes \CI(B)$ is dense in $\Psi^{m}(\GR \times B)/
\Psi^{m-1}(\GR \times B)$ in the corresponding Frechet topology
(defined by the isomorphism of Theorem \ref{Theorem.Onto}).
It follows that $\Psi^{m}(\GR \times B)$ consists of smooth families of
operators in $\Psi^{m}(\GR)$ parametrized by $B$, see \cite{Atiyah-Singer4},
page 122 and after, where families of pseudodifferential operators are 
discussed.

We obtain that $A(\GR \times B)$ is to pull back of  $A$ to $\Gr0 \times B$.
\end{example}

The following example generalizes the tangent groupoid of Connes; here
we  closely follow \cite[II,5]{Connes3}. The groupoid defined
below also appears in \cite{we:blowing} and is related to the notion
of explosion of manifolds.

\begin{example} The adiabatic groupoid $\GR_{\adb}$ 
associated to $\GR$ is defined as follows. The space 
of units is $\GR_{\adb}^{(0)}=[0,\infty) \times \Gr0$ with the product
manifold structure. The set of arrows $\GR_{\adb}^{(1)}$ is 
defined to be the disjoint union $A(\GR) \cup  (0,\infty) \times \Gr1$,
and $d(t,g)=(t,d(g))$, $r(t,g)=(t,r(g))$ if $t>0$,
$d(v)=r(v)=(0,x)$ if $v \in T_x\GR_x$. 
The composition is $\mu(\gamma,\gamma')=(t,gg')$
if $\gamma=(t,g)$ and $\gamma'=(t,g')$ for $t>0$
(necessarily the same $t$!) and $\mu(v,v')=v +v'$
if $v,v' \in T_x\GR_x$. 

The smooth structure on the set of arrows is the product structure for $t > 0$.
In order to define a coordinate chart at a point $v \in T_x\GR_x$
choose first a coordinate system $\psi:U=U_1 \times U_2 \to \Gr1$, 
$U_1  \subset \RR^p$ and $U_2\subset \RR^n$ being open sets
containing the origin, $U_2$ convex, with the following properties: 
$\psi(0,0)=x \in \Gr0 \subset \Gr1$, $d(\psi(s,y_1))=d(\psi(s,y_2))=\phi(s)$
and $\psi(U) \cap \Gr0=\psi(U_1 \times \{0\})$. Here
$\phi:U_1 \to \Gr0$ is a coordinate chart of $x$ in $\Gr0$.
We identify, using the differential
$\operatorname{D}_2\! \psi$ 
of the map $\psi$, the vector space $\{s\} \times \RR^n$
and the  tangent space $T_{\phi(s)}\GR_{\phi(s)}=A_{\phi(s)}(\GR)$. 
We obtain then  coordinate charts 
$\psi_\epsilon:[0,\epsilon) \times U_1 \times \epsilon^{-1}U_2 
\to  \Gr1$, $\psi_\epsilon(0,s,y)=(0,(\operatorname{D}_2\!\psi) (s,y)) \in 
T_{\phi(s)}\GR_{\phi(s)}=A_{\phi(s)}(\GR)$ and 
$\psi_\epsilon(t,s,y)=(t,\psi(s,ty)) \in (0,1) \times \Gr1$. For 
$\epsilon$ very small the range of $\psi_\epsilon$ will contain $v$.

For $\GR=M\times M$ as in the first example the groupoid
$\GR_{\adb}$ is the tangent groupoid defined by  Connes, and the
algebra of pseudodifferential operators is the algebra of asymptotic 
pseudodifferential operators \cite{Widom1}.
In general an operator $P$ in $\Psi^m(\GR_{\adb})$
will restrict to an adiabatic  family  $P=(P_{t,x}, t>0, x \in \Gr0)$
which will have an ``adiabatic limit'' at $t=0$ given by the operator
$P$ at $t=0$. 

The Lie algebroid of $\GR_{\adb}$ is the adiabatic Lie algebroid
associated to $A(\GR)$, $ A(\GR_{\adb})=A(\GR)_{t}$ (using the notation
of Theorem \ref{Theorem.Quantization}).
This gives a procedure for integrating adiabatic Lie algebroids.
Using pseudodifferential operators on the adiabatic groupoid we obtain 
an explicit quantization of symbols on $A^*$
generalizing  Theorem \ref{Theorem.Quantization}. The proof proceeds 
exactly in the same way.
\end{example}

\begin{theorem}
\label{Theorem.Quantization2}
The inverse limit $\projlim \Psi^\infty(\GR_{\adb})/t^n\Psi^\infty(\GR_{\adb})$
is a deformation  quantization of the commutative algebra 
${\mathcal S}_{\cl}^\infty(A^*(\GR))$ of classical symbols.
\end{theorem}

The space ${\mathcal S}_{\cl}^\infty(A^*(\GR))$ appearing in the statement of
the theorem above is the union of all symbol spaces ${\mathcal S}_{\cl}^m(A^*(\GR))$
and is a commutative algebra under pointwise multiplication.

Of course Theorem \ref{Theorem.Quantization} provides us with a
$*$-product whose multiplication is given by differential operators,
and hence this $*$-product extends to all smooth functions (even to
functions defined on open subsets). The usefulness of the theorem
above is that it gives in principle a {\em nonperturbative} (i.e. not
just formal) deformation quantization, close in
spirit to that of {\em strict} deformation quantization introduced by
Rieffel \cite{Rieffel1}.  The algebras $A_t$ \cite{Rieffel1} turn out
to be groupoid algebras.

\begin{example} 
This example provides a treatment in our settings of 
the $b$- and $c$-calculi defined by Melrose 
\cite{Melrose42,Melrose46,Melrose-Nistor2,Melrose-Nistor1}
on a manifold with boundary $M$.
 
Define first a groupoid $\GR_\phi(M)$ associated to $M$  and
an increasing diffeomorphism $\phi:\RR \to (0,\infty)$ as follows.
If $M=[0,\infty)$ the action by translation of $\RR$ on itself extends 
to an action on $M$ fixing $0$, not smooth in general, 
defined using the isomorphism $\phi$. 
Define $\GR_\phi(M)$ to be the transformation groupoid associated this
action of $\RR$ on $M$. If $M=[0,1)$
then $\GR_\phi(M)$ is defined to be the reduction $\GR_\phi(M) = 
\GR_\phi([0,\infty)) \cap d^{-1}(M) \cap r^{-1}(M)$ of 
$\GR_\phi([0,\infty))$ to $[0,1)$.

Suppose next that $M = \pa M \times [0,1)$. We then define
$\GR_\phi(M) = \GR_\phi([0,1)) \times (\pa M \times \pa M)$
where $\pa M \times \pa M$ is the equivalence groupoid of $\pa M$ considered
in the first example. 
For an arbitrary manifold with boundary $M$ write 
$M = M_0 \cup U$ where $U = M \setminus \pa M$ 
and $M_0$ is diffeomorphic to $\pa M \times [0,1)$. (Our construction
will depend on this diffeomorphism.)
Then we define $\Gr0_\phi(M)=M$ and 
$\Gr1_\phi(M) = \Gr1_\phi(M_0) \cup (U \times U)$ with the induced operations. 

If $\phi(t)=e^t$ or $\phi(t)=-t^{-1}$ (for $t<<0$) then $\GR=\GR_\phi$
will be an almost differentiable groupoid and we obtain
$\tPS{m} \subset \Psi_b(M)$ in the first case and  
$\tPS{m} \subset \Psi_c(M)$ in the second case. The first groupoid
does not depend on any choices.
\end{example}

\section{Distribution kernels}\label{D.k}

In this section we characterize the reduced (or convolution)
distribution kernels of operators in $\PS{m}$ following  \cite{Melrose42} 
(see also \cite{Hormander3}) as compactly supported
distributions on $\GR$, conormal to the set of units $\Gr0$. 

Denote by $\ENDG$ the bundle $\Hom(d^*(E),r^*(E))=r^*(E)\otimes d^*(E)'$
on $\Gr1$,
where
$V'$ denotes as usual the dual of the vector bundle $V$. Using the
relations $d \circ \iota =r$ and $r \circ \iota =d$ we see that
$\ENDG$ satisfies
\begin{equation}\label{eq.selfdual}
\iota^*(\ENDG) \simeq d^*(E) \otimes r^*(E)' \simeq \ENDG'.
\end{equation}

We define  a convolution product on the space 
$\CIc(\Gr1, \ENDG \otimes d^*({\mathcal D}))$ of compactly supported
smooth sections of the bundle $\ENDG \otimes d^*({\mathcal D}))$ by the formula
\begin{equation}
f_1* f_2(g) = \int_{\{(h_1,h_2),h_1h_2=g\}} f_1(h_1) f_2(h_2)\;. 
\end{equation}
The multiplication on the right hand side is the composition of 
homomorphisms giving a linear map
\begin{gather*}
\Hom(E_{d(h_1)},E_{r(h_1)}) \otimes \Hom(E_{d(h_2)},E_{r(h_2)})
\otimes \mathcal{D}_{d(h_1)}\otimes \mathcal{D}_{d(h_2)}
\longrightarrow \\
\Hom(E_{d(g)},E_{r(g)}) \otimes \mathcal{D}_{d(h_1)}\otimes 
\mathcal{D}_{d(h_2)}\\ 
f_1(h_1) \otimes f_2(h_2) \longrightarrow f_1(h_1) f_2(h_2) ,
\end{gather*}
defined since $d(h_1)=r(h_2)$.
To see that the integration is defined we parametrize the set 
$\{(h_1,h_2),h_1h_2=g\}$ as $\{(gh^{-1},h),h \in \GR_{d(g)}\}$ which shows
that this set is a smooth manifold, 
and notice that we can invariantly define the integration with respect to 
$h$ taking advantage of the $1$-density factor 
${\mathcal D}_{d(h_1)}={\mathcal D}_{r(h)}=(\Omega_d)_{h}$.
If we choose a hermitian metric on ${\mathcal D}^{-1/2} \otimes E$,
we obtain a conjugate--linear involution 
(making $\CIc(\Gr1, \ENDG \otimes d^*({\mathcal D}))$ into
a $*$-algebra). 

Consider an operator $P=(P_x,x \in \Gr0) \in \PS{-\infty}$ and let 
$k_x$ be the distribution kernel of $P_x$, a smooth section 
$k_x \in \CI(\GR_x \times \GR_x;r_1^*(E) \otimes r_2^*(E)' \otimes \Omega_2)$, 
using the notation $\Omega_2=p_2^*(\Omega_d)=r_2^*(\mathcal{D})$ of \eqref{kernel}.
We define the reduced distribution kernel $k_P$ of the smoothing operator $P$ by 
\begin{equation}
\label{def.smooth.reduced.kernel}
k_P(g)=k_{d(g)}(g, d(g)) \in E_{r(g)} \otimes E_{d(g)}' 
\otimes \mathcal{D}_{d(g)}\,.
\end{equation}
This definition will be later extended to all of $\PS{\infty}$.

The following theorem is one of the main reasons we consider
{\em uniformly} supported operators.

\begin{theorem}
\label{Theorem.Residual.Isom} 
The reduced kernel map $P \to k_P$ \eqref{def.smooth.reduced.kernel} 
defines an isomorphism of the residual ideal $\PS{-\infty}$ with the 
convolution algebra $\CIc(\Gr1, \ENDG \otimes d^*({\mathcal D}))$.
\end{theorem}

\begin{proof} Let $P$ and $k_x$ be as above. We know from Lemma
\ref{Lemma.Smoothness} that the collection of all sections 
$k_x$ defines a {\em smooth}
section of $r_1^*(E)\times r_2^*(E)'\otimes \Omega_2$
over the manifold $\{(g_1,g_2),d(g_1)=d(g_2)\}$.
The relation $P_{r(g)}U_g = U_g P_{d(g)}$ gives the invariance relation
$k_{r(g)}(h',h)=k_{d(g)}(h'g,hg) \in
E_{r(h')} \otimes E_{r(h)}' \otimes \mathcal{D}_{r(h)}$ 
for all arrows $g \in \Gr1$,
and $h,h' \in \GR_{r(g)}$.
It follows that $k_{d(h)}(h',h)=k_{r(h)}(h'h^{-1},r(h))
=k_P(h'h^{-1})$. The section $k_P$ is 
well defined, smooth and completely  determines all kernels $k_x$ and hence also 
the operator $P$. Moreover the section $k_P$ has compact support
because $\supp(k_P)=\supp_\mu(P)=
\mu\circ(id \times \iota)\big(\overline{\cup_x \supp(k_{x})}\big)$ and 
the reduced support $\supp_\mu(P)$ 
of $P$ is compact since $P$  is uniformly 
supported. The distribution kernel $k^{PQ}_x$ of the product
$P_xQ_x$ of two operators $P_x,Q_x \in \Psi^{-\infty}(\GR_x)$ is 
$$
k^{PQ}_x(g,g'')=
\int_{\GR_x}k^P_x(g,g')k^Q_x(g',g'')dg'
$$
where $k^P_x$ and $k^Q_x$ are the distribution kernels
of $P_x$ and, respectively, $Q_x$. From this, taking into account the
definitions of $k_{PQ}$, $k_P$ and $k_Q$, we obtain 
$$
k_{PQ}(g)=
\int_{\GR_x}k_P(g{g'}^{-1})k_Q(g')dg'\ .
$$
This means that $k_{PQ}=k_P * k_Q$ and hence the  map $P \to k_P$ 
establishes the desired
isomorphism.
\end{proof}

We will now use duality  to extend the definition of the reduced
distribution 
kernel to any operator $P \in \PS{\infty}$. 
Let ${\mathcal L}=\Omega_{\Gr0}$ be the line bundle of $1$-densities 
on $\Gr0$ and
$\VD=\Omega_d\vert_{\Gr0}$ be the bundle of vertical $1$-densities as above.
Define
\begin{gather*}
\PS{-\infty}_{\mathcal L}=\PS{-\infty}\otimes_{\CI(\Gr0)} 
\CI(\Gr0,{\mathcal D}^{-1} \otimes {\mathcal L})\simeq \\
\PS{-\infty} \otimes_{\CI(\Gr0)} \CI(\Gr0,{\mathcal D}^{-1}) 
\otimes_{\CI(\Gr0)}\CI(\Gr0,{\mathcal L})
\end{gather*}
where the tensor products are defined using
the inclusion $\CI(\Gr0)\subset\PS\infty $. We note that the bundle
${\mathcal L}$ plays an important role in connection with  the modular class of
a groupoid  \cite{ELW}, since it carries a natural representation of
the groupoid.

The relation $k_{Pf}(g)=k_P(g)f(d(g))$ for $f \in \CI(\Gr0)$ and 
$P \in \PS{-\infty}$   give using Theorem \ref{Theorem.Residual.Isom} the 
isomorphism
\begin{equation}\label{eq.equation.L}
\PS{-\infty}_{\mathcal L} \simeq \CIc(\Gr1, 
\ENDG \otimes d^*({\mathcal L})).
\end{equation}
The space $\PS{-\infty}_{\mathcal L}$ comes equipped with a 
natural linear functional $\TR$ such that if
$P_0\in \PS{-\infty}$, $\xi \in \CI(\Gr0,{\mathcal D}^{-1})$
and $\nu \in \CI(\Gr0,{\mathcal L})$ then 
\begin{equation*}
\TR(P_0 \otimes\xi \otimes \nu)=\int_{\Gr0}tr(k_{P_0}(x)\xi(x)) d\nu(x)
\end{equation*}
defined by integrating the function $tr(k_{P_0}(x)\xi(x))$ with respect to 
the $1$-density (\ie measure) $\nu$. An operator 
$P \in \PS{m}$ defines a continuous linear functional (\ie distribution)
$k_P^{\iota}:\PS{-\infty}_{\mathcal L} \to \CC$
by the formula
$k_P^{\iota}(P_0 \otimes\xi \otimes \nu)=\TR(PP_0\otimes\xi \otimes \nu)$.
It is easy to see using Equation \eqref{eq.selfdual} that the map 
$f \to \tilde f = f \circ \iota$, $\iota(g)=g^{-1}$, establishes isomorphisms
\begin{gather}\label{eq.distr.isom}
\Phi:\CIc(\Gr1, \ENDG \otimes d^*({\mathcal L})) 
\stackrel{\iota^*}{\longrightarrow} 
\CIc(\Gr1, \ENDG'\otimes r^*({\mathcal L})) \simeq \\ \nonumber
\CIc(\Gr1, (\ENDG \otimes d^*({\mathcal L}))' \otimes d^*({\mathcal L})
\otimes r^*({\mathcal L}))
\simeq \\ \nonumber
\CIc(\Gr1, (\ENDG \otimes d^*({\mathcal L}))' \otimes \Omega_\GR)
\end{gather}
whose composition we denote by $\Phi$, so that
$\Phi(P_0\otimes\xi \otimes \nu)=(k_{P_0}\xi\nu)\circ \iota = 
\iota^*(k_{P_0}\xi\nu)$. We obtain in this way from 
$k_P^\iota$ a distribution 
$k_P \in {\mathcal C}^{-\infty}(\Gr1;\ENDG \otimes d^*({\mathcal D}))$
defined by the formula 
\begin{equation}\label{eq.def.kp}
\langle k_P, f \rangle = k_{P}^\iota (\Phi^{-1}(f)).
\end{equation}
An other way of writing the formula above is
\begin{equation}\label{eq.def.kp2}
\langle k_P, \iota^*(k_{P_0}\xi\nu) \rangle = \TR(PP_0\otimes\xi \otimes \nu)=
\int_{\Gr0}tr(k_{PP_0}(x)\xi(x)) d\nu(x).
\end{equation}

\begin{proposition} If $P\in \PS{-\infty}$ is a regularizing 
operator then the kernels $k_P$ defined in Equations 
\eqref{def.smooth.reduced.kernel} and \eqref{eq.def.kp} coincide.
\end{proposition}

\begin{proof} To make a distinction for the purpose of this proof,   denote by
$k_P^{\dist}$ the distribution defined by \eqref{eq.def.kp}.
Let $\nu$ be a smooth section of ${\mathcal L}$,
$\xi$  a smooth section of ${\mathcal D}^{-1}$
and $P,P_0 \in \PS{-\infty}$. Using Equation \eqref{eq.def.kp2}
we obtain
\begin{multline*}
\langle k_P^{\dist}, \iota^*(k_{P_0}\xi\nu) \rangle = 
\TR(PP_0\otimes\xi \otimes \nu)=
\int_{\Gr0}tr(k_{PP_0}(x)\xi(x)) d\nu(x)=\\
\int_{\Gr0} \left(\int_{\GR_{x}} 
tr\big(k_{P}(h^{-1}) k_{P_0}(h)\xi(x)\big) \right)d\nu(x)
=\langle k_P \circ \iota,k_{P_0} \xi \nu\rangle
=\langle k_P ,\iota^*(k_{P_0} \xi \nu)\rangle.
\end{multline*}
\end{proof}

\begin{definition} The distribution $k_P\in 
{\mathcal C}^{-\infty}(\Gr1;\ENDG \otimes d^*({\mathcal D}))$, defined for any
operator $P \in \PS{m}$ by Equation \eqref{eq.def.kp}
will be called the reduced (or convolution) distribution 
kernel of $P$, or simply the reduced kernel of $P$, and will be denoted $k_P$.
\end{definition}

We now relate the action of $\PS{\infty}$ by multiplication on 
$\PS{-\infty}$, respectively on $\PS{-\infty}_{\mathcal L}$,
to that on $\CIc(\GR,r(E))$.
\begin{gather}\label{eq.isomorphisms1}
\PS{-\infty} \simeq \CIc(\Gr1; r^*(E))
\otimes_{\CI(\Gr0)} \Gamma(E') 
\otimes_{\CI(\Gr0)} \Gamma({\mathcal D})\\ 
\label{eq.isomorphisms2}
\PS{-\infty}_{\mathcal L} \simeq \CIc(\Gr1; r^*(E))
\otimes_{\CI(\Gr0)} \Gamma(E') 
\otimes_{\CI(\Gr0)} \Gamma({\mathcal L})
\end{gather}
such that the left action by multiplication of
$\PS{\infty}$ on $\PS{-\infty}$ becomes
$P(f\otimes \eta \otimes \xi)=Pf \otimes \eta \otimes \xi$
where $\eta$ is a smooth section of $E'$ and $\xi$
is a smooth section of ${\mathcal D}$ or ${\mathcal L}$. Moreover
the kernel of $P_0=f\otimes \eta \otimes \xi$ is 
$k_{P_0}(g)=f(g)\otimes \eta(d(g))\xi(d(g))$. Thus in order to 
define the distribution $k_P$, for arbitrary $P$, it is enough
to compute $\TR(Pf_0\otimes\eta\otimes\nu)$ where $\nu$ is a 
density.

Fix a unit $x$ and choose a coordinate chart $\phi:U_0 \to U \subset \Gr0$ where 
$U_0$ is an open subset of $\RR^k$ containing $0$, 
$k =\dim \Gr0$ and $\phi(0)=x$. By decreasing
$U_0$ if necessary we can assume that the tangent space $T\Gr0$ is
trivialized over $U$.
Consider the  diffeomorphism $\exp_\nabla:V_0 \to V \subset \Gr1$ associated
to a right invariant connection $\nabla$ as in 
\eqref{eq.invariant.connection} and \eqref{eq.exp.diffeomorphism}
where $V_0 \subset A(\GR)$ is an open neighborhood of the zero section. 
It maps the zero section of $A(\GR)$ to $\Gr0$. Choose a connection 
on $E$ which lifts to an invariant connection on $r^*(E)$.
By decreasing $V$ if necessary and using  the invariant connection on $r^*(E)$
we obtain canonical trivializations of  $r^*(E)$ on each fiber
$V \cap \GR_x$.
Denote by $\theta_h:E_{r(h)}\otimes E_{d(h)}' \to \End(E_x)$ the isomorphism 
induced by the connection $\nabla'$ (defined using parallel transport along
the geodesics of $\nabla$) where $x=d(h)$ and $h$ is in $V$. Decreasing further
$V$ and $U$ we can assume that $E$ is trivialized over $d^{-1}(U) \cap V$
and that $\phi$ and $\exp_\nabla$ give a fiber preserving diffeomorphism 
$\psi:U_0 \times W \to d^{-1}(U) \cap V$ where $W\subset \RR^n$ 
is an open set, identified with an open neighborhood  of the zero section 
in $T_x\Gr0$. The diffeomorphism $\psi$ we have just constructed satisfies 
$d(\psi(s,y))=\phi(s)$. The maps $\psi$ and $\theta_h$ yield isomorphisms 
\begin{gather}\label{distr.isom}
{\mathcal C}^{-\infty}(d^{-1}(U) \cap V, \ENDG \otimes d^*({\mathcal D})) 
\simeq \\ \nonumber {\mathcal C}^{-\infty}
(U_0\times W, \psi^*(\ENDG \otimes d^*({\mathcal D})))
\simeq {\mathcal C}^{-\infty}(U_0\times W,E_x \otimes E_x')
\end{gather}
whose composition is denoted $\Theta_\psi$.

Next theorem describes the reduced distribution kernels $k_P$ of operators
$P$ in  $\PS{m}$. We use the notation introduced above.

\begin{theorem}\label{Theorem.Full.Isom} 
For any operator $\FAM \in\PS{m}$ the reduced distribution kernel $k_P$ 
satisfies: 

(i) If $\psi:U_0 \times W \to V_1\subset V$, $W\subset \RR^n$ open, 
is a diffeomorphism satisfying $\psi(s,0)=d(\psi(s,y))$,
then there exists a symbol 
$a_P\in {\mathcal S}_{\cl}^m(U_0 \times \RR^n;\End(E_x))$, 
such that $k_P\circ \iota=\Theta_\psi^{-1}(k)$
on $V_1$, where $k$ is the distribution
\begin{equation*}
k(s,y)=(2\pi)^{-n}\int_{\RR^n} e^{-iy\cdot \zeta}a_P(s,\zeta)d \zeta \in  \End(E_x)\,, 
\end{equation*}
the integral being an oscillatory integral. Moreover, after
suitable identifications, $a_P$ is a representative of the
principal symbol of $P$.

(ii) The singular support of $k_P$ is contained in $\Gr0$. 

(iii) The support of $k_P$ is compact, more precisely $\supp(k_P)=\supp_\mu(P)$.

(iv) For every distribution $k\in {\mathcal C}^{-\infty}(\GR;E_0)$, 
satisfying the three conditions above, there exists 
$P \in \PS{m}$ such that $k=k_P$.
\end{theorem}

Note that $k_P\circ \iota$, $\iota^*(k_P)$ and $k_P^{\iota}$ all
denote the same distribution.

\begin{proof} Write $\phi(s)$ for $\psi(s,0)=d(\psi(s,y)).$
According to Definitions 
\ref{Def.Differentiable} and \ref{Main.definition}, there exists  
a classical symbol $a \in {\mathcal S}_{\cl}^m(U_0 \times T^*W;\End(E_x))$
such that $P_{\phi(s)}=a(s,y,D_y)$ (modulo regularizing operators)
on $\GR_{\phi(s)} \cap V_1\simeq W$. 

Let $P_0\in \PS{-\infty}$, $\xi \in \CI(\Gr0,{\mathcal D}^{-1})$
and $\nu \in \CI(\Gr0,{\mathcal L})$. Assume, using the isomorphisms
\eqref{eq.isomorphisms1} and \eqref{eq.isomorphisms2}, that 
${P_0}=f_0 \otimes \eta$ where $\eta \in \Gamma(E')= \CI(\Gr0,E')$ 
and $f_0\in \CI(\GR,r^*(E)\otimes \Omega_d)$,
so that $f_0\xi$ is a section of $\CIc(\GR,r^*(E))$. Then we have
$$
tr(k_{PP_0}(x)\xi(x))=
\eta(P_x(f_0\xi\vert_{\GR_x})(x)).
$$
Suppose $f_0$ is supported in $V_1$,
and denote by $f_s$ the section of $E_x$ that corresponds to 
$f_0\xi\vert_{\GR_{\phi(s)}}$ under the diffeomorphism 
$\GR_{\phi(s)} \cap V_1 \simeq \{s\} \times W = W$ induced by $\psi$.
We then have
\begin{gather*}
\langle \iota^*(k_P), k_{P_0}\xi\nu \rangle =
\langle k_P, \iota^*(k_{P_0}\xi\nu) \rangle = 
\TR(PP_0\otimes\xi \otimes \nu)=\\
\int_{\Gr0}tr(k_{PP_0}(x)\xi(x)) d\nu(x)=
\int_{\Gr0\cap U}\eta\big(P_x(f_0\xi\vert_{\GR_x})(x)\big)d\nu(x)=\\
\int_{U_0}\eta\left(\int_{\RR^n}\int_{W}
e^{-i y\cdot \zeta}a(s,0,\zeta)f_s(y)dy
d\zeta \right)d\nu(s)=\\
\int_{U_0}\int_{\RR^n}\int_W  e^{-i y\cdot \zeta}
tr\big(a(s,0,\zeta)f_s(y)\otimes \eta \big)dyd\zeta d\nu(s)=\\
\int_{U_0\times W}  tr\left(f_s(y)\otimes \eta\int_{\RR^n}e^{-i y\cdot \zeta}
a(s,0,\zeta)d\zeta  \right)dyd\nu(s)
\end{gather*}
where the first integral is really a pairing between the distribution $k$
obtained from (i) for $a_P(s,\zeta)=a(s,0,\zeta)$,
and the smooth section $f_s \otimes \eta \otimes \nu$.
Since $\End(E_x)$ is canonically its own dual this shows that the 
distribution $k_P$ is the conormal distribution to
$\Gr0$ given by (i). 

To prove (iii) and (iv) observe that $k_x$ is the restriction to
$\GR_x \times \GR_x$ of the distribution $\mu_1^*(k_P)$, where
$\mu_1(h',h)=h'h^{-1}$ and the distribution $\mu_1^*(k_P)$ is
defined by $\langle \mu_1^*(k_P),f\rangle 
=\langle k_P(g) ,\int_{h_1h_2=g} f(h_1,h_2)\rangle.$ Then we can
define $P_x$ by its distribution kernel $k_x$. From (i) it follows that
$k_x$ is conormal to the diagonal and hence $P_x$ is a pseudodifferential
operator. 

In order to check (ii) fix $g \not \in \Gr0$ and let $\varphi$ be a smooth cut-off
function, $\varphi=1$ in a neighborhood of $\Gr0$, $\varphi=0$ in a neighborhood
of $g$. Consider again the distribution 
$\mu_1^*((1-\varphi)k_P)=(1-\varphi\circ \mu_1)\mu_1^*(k_P)$. Its restriction
to $\GR_x$ is $(1-\varphi\circ \mu_1)k_x$ 
which is smooth since the singular
support of $k_x$ ($=$ the distribution kernel of $P_x$) is contained in the
diagonal of $\GR_x \times \GR_x$, and $1-\varphi\circ \mu_1$ vanishes there.
It follows that $\mu_1^*((1-\varphi)k_P)$ is smooth 
and hence $(1-\varphi)k_P$ is also smooth.
\end{proof}

\begin{corollary}The distribution $k_P$ is conormal at 
$\Gr0$ and smooth everywhere else. In particular
the wave-front set of $k_P$ is a subset of the annihilator of $T\Gr0$:
$WF(k_P)\subset (T\GR/T\Gr0)^*\subset T^*\GR\vert_{\Gr0}$.
\end{corollary}

\begin{proof}This is a standard consequence of (i) and (ii) 
above, see \cite{Hormander3}, section 12.2.
\end{proof}

We remark that $(T\GR/T\Gr0)^*$ is naturally identified with $A^*(\GR)$.
Denote by ${\mathcal S}_{c}^m(A^*(\GR);\End(E))\subset
{\mathcal S}_{\cl}^m(A^*(\GR);\End(E))$ the space of classical symbols
with support in a set of the form $\pi^{-1}(K)$ where $\pi:A^*(\GR) \to
\Gr0$ is the projection and $K \subset \Gr0$ is a compact subset.

\begin{corollary}\label{cor.support} 
Let $V$ be a neighborhood of $\Gr0$ in $\GR$. Then any
$P \in \PS{m}$ can be written as $P=P_1 +P_2$ where $P_1$ has reduced support 
$\supp_\mu(P_1)$ contained in $V$ and $P_2 \in \PS{-\infty}$.
\end{corollary}

\begin{proof}Let $\phi$ be a smooth cut-off function, equal
to $1$ in a neighborhood of $\Gr0$ and with support in  $V$.
Define $P_2 \in \PS{-\infty}$ by  $k_{P_2}=k_P(1-\phi)$.
This is possible using Theorem \ref{Theorem.Residual.Isom}
since by (ii) of the theorem above $k_P(1-\phi)$
is a smooth compactly supported section of an appropriate bundle. 
Then $P_1=P-P_2$ and $P_2$ satisfy the requirements of the statement.
\end{proof}

\begin{theorem}\label{Theorem.Onto}
The principal symbol map $\sigma_m$ in Equation \eqref{Eq.Def.PS} is onto; 
hence it establishes an isomorphism
$$\PS{m}/\PS{m-1} \simeq 
{\mathcal S}_{c}^m(A^*(\GR);\End(E))
/{\mathcal S}_{c}^{m-1}(A^*(\GR);\End(E))$$
for any $m$.
\end{theorem}

\begin{proof} We only need to prove that $\sigma_m$
is onto. If follows from the proof of Theorem \ref{Theorem.Full.Isom}
that $\sigma_m(P)$ is the class of the symbol $a_P$ appearing
in equation in (i). Given a symbol 
$a \in {\mathcal S}_{c}^m(A^*(\GR);\End(E))$
the equation in (i) defines a distribution $k_0$ in a small neighborhood
of $\Gr0$ in $\GR$. Using a smooth cut-off function we obtain
a distribution $k$ on $\GR$ that coincides 
with $k$ in a neighborhood of $\Gr0$ and is smooth outside
$\Gr0$. From (iv) we conclude that there exists an
operator $P$ with $k_P=k$ which will then 
necessarily satisfy $\sigma_m(P)=a+{\mathcal S}_{c}^{m-1}(A^*(\GR);\End(E))$.
\end{proof}

\section{The action on sections of $E$}

In this section we define a natural action of $\PS{m}$ on 
sections of $\Gr0$, thus generalizing the action of classical
pseudodifferential operators on functions.

Let $\phi$ be a smooth section of $E$ over $\Gr0$. Define
\begin{equation}
\label{tilde1}
\tilde \phi \in \CI(\Gr1,r^*(E))\,,\,
\tilde\phi(g)=\phi(r(g)).
\end{equation} 

\begin{lemma} 
If $\FAM$ belongs to $\PS{\infty}$, then for any 
section $\phi$ in $\CI(\Gr0,E)$ there exists a unique section
$\psi \in \CI(\Gr0,E)$ such that  $P\tilde \phi = \tilde \psi$.
\end{lemma}

\begin{proof}
Observe first that given a section $\gamma$ of $r^*(E)$ 
over $\Gr1$ we can find a section 
$\phi$ of $E$ over $\Gr0$ such that $f=\tilde \phi$ if and only if 
$f(g'g)=f(g')$ for all $g$ and $g'$, i.e. if and only if
\begin{equation}
\label{tilde}
U_g f_x =f_y\, , \, \text{ for all } g,\,x,\,y \text{ such that }x=d(g)
\text{ and } y=r(g).
\end{equation}
We then have 
$$U_g \tilde \phi_x = \tilde\phi_y \Rightarrow 
P_y U_g \tilde \phi_x = P_y \tilde\phi_y \Rightarrow 
U_g P_x \tilde \phi_x = P_y \tilde\phi_y \Rightarrow 
U_g (P \tilde \phi)_x = (P \tilde\phi)_y$$
and hence $P\tilde \phi$ satisfies  \eqref{tilde}. Thus
we can find a section $\psi$ of $E$ over $\Gr0$ such that
$P\tilde \phi = \tilde \psi$. Note that $P_x \tilde \phi_x$ is defined
since $P_x$ is properly supported.

The uniqueness of the section $\psi$ follows from the fact that
the map $\phi \to \tilde \phi$ is one-to-one, and the smoothness of
$\psi$ follows from Lemma \ref{Lemma.Differentiable}.
\end{proof}

The representation given by the following theorem reduces to
the trivial representation in the case a group (see also comments bellow).

\begin{theorem} 
There exists a canonical representation $\pi_0$  of the algebra
$\PS{\infty}$ on $\CI(\Gr0,E)$ given by $\pi_0(P)\phi =\psi$
where, using the notation of the previous lemma,
$\psi$ is the unique section satisfying $\tilde \psi=P \tilde \phi$. 
Moreover $\pi_0(P)$ maps compactly supported sections to compactly supported
sections.
\end{theorem}

\begin{proof} The fact that $\pi_0$ is well defined follows from 
the uniqueness part of the previous lemma. It is clearly a
representation. We only need to check that $\pi_0(P)$ maps compactly
supported sections to compactly supported sections. Let $L_1 \subset
\Gr0$ be the support of $\phi$, $L_2=\supp (P)$. Then the support of
$\pi_0(P)$ is contained in $L_2 L_1$.
\end{proof}

Assume that $E$ is a trivial line bundle. Then $\CI(M)$ and $\Gamma(A)$
act naturally on $\CI(M)$ and this action  satisfies the relations
\eqref{eq.relations} which means that it gives rise to a representation
of $U(A)=\operatorname{Diff}(\GR)$ on $\CI(M)$. Then $\pi_0$
is an extension of this representation. If $\GR=G$ is a group, then
$\pi_0$ extends the trivial representation.  In order to generalize this fact to
arbitrary representation of $G$ we need the following definition.
 
\begin{definition}
\label{equivariant}
An equivariant bundle $(V,\rho)$ on $\Gr0$ is a differentiable vector bundle 
$E$ together with a bundle isomorphism $\rho: d^*(V) \longrightarrow r^*(V)$ 
satisfying $\rho(gh)=\rho(g)\rho(h)$.
\end{definition}

An equivariant bundle is also called a representation of $\GR$.
Given an equivariant bundle $(V,\rho)$, we can define a representation 
$\pi_\rho$ of
the groupoid algebra $\CIc(\GR, d^*({\mathcal D}))$ on $\CIc(\Gr0,V)$
by the formula
\begin{equation}
\label{representation}
(\pi_\rho(f)\phi)(x)=\int_{\GR_x} f(h^{-1})\rho(h^{-1})\phi(r(h)).
\end{equation}
Note that the integration is defined and gives an element
of $V_x$ since $f(h^{-1})\phi(r(h))$
is in $\CIc(r^*(V) \otimes \Omega_d)$ and hence 
that $f(h^{-1})\rho(h^{-1})\phi(r(h))$ is a smooth compactly supported section
of $d^*(V) \otimes \Omega_d$.

The following proposition has no obvious analog in the classical
theory because the pair groupoid has no nontrivial representations. If
one moves one step up and considers the fundamental groupoid,
nontrivial representations exist, and the following lemma says that
geometric operators (\ie the ones that lift to the universal covering
space) act on sections of flat bundles. A representation of a groupoid
thus resembles a flat bundle.

\begin{proposition}Let $(V,\rho)$ be an equivariant bundle and $E$
an arbitrary bundle on $\Gr0$. There exists a natural morphism 
$T_\rho:\PS{\infty} \to \Psi^{\infty}(\GR;V\otimes E)$
and hence there exist a canonical action $\pi_\rho=\pi_0 \circ T_\rho$ 
of $\tPS{\infty}$ on  $\CI(\Gr0,E\otimes V)$ and $\CIc(\Gr0,E\otimes V)$  
which extends the representation defined in \eqref{representation}.
\end{proposition}

\begin{proof} Let 
$$W_{\rho,x}:\CIc(\GR_x;r^*(E)) \otimes V_x =\CIc(\GR_x;r^*(E) \otimes
d^*(V)) \to 
\CIc(\GR_x;r^*(E \otimes V))$$ be the isomorphism defined by $\rho$ as in the
definition \ref{equivariant}. It is easy to see that this gives an isomorphism
$W_\rho:\CIc(\GR;r^*(E) \otimes d^*(V)) \to
\CIc(\GR;r^*(E \otimes V))$. Define an operator on $\CIc(\GR_x;r^*(E \otimes V))$
by the formula
$$(T_\rho(P))_x = W_{\rho,x}(P_x \otimes id_{V_x})W_{\rho,x}^{-1}.$$
The relation $W_{\rho,x} (U_g \otimes \rho(g))=U_g W_{\rho,x}$
shows that the family $(T_\rho(P))_x$, $x \in \Gr0$ satisfies the invariance
condition $(T_\rho(P))_xU_g=U_g(T_\rho(P))_y$, for $d(g)=x$ and $r(g)=y$.
The uniform support condition is satisfied since $\supp(T_\rho(P))=\supp_\mu(P)$.
It follows that the family $(T_\rho(P))_x$ defines an operator $T_\rho(P)$
in $\Psi^{\infty}(\GR;V\otimes E)$. The multiplicativity condition 
$T_\rho(PQ)=T_\rho(P)T_\rho(Q)$ follows from definition and hence $T_\rho$
is a morphism.
\end{proof}


\end{document}